\documentclass[titlepage,11pt]{article}
\usepackage{latexsym}
\usepackage{amsfonts}
\usepackage{amsmath}
\usepackage{tikz}
\usetikzlibrary{decorations.pathreplacing,calligraphy}
\newcommand\blackslug{\hbox{\hskip 1pt \vrule width 4pt height 8pt depth 1.5pt
        \hskip 1pt}}
\newcommand\bbox{\hfill \quad \blackslug \bigbreak}

\def\DD{\hbox{-}}
\def\LL{,\ldots,}
\def\cupcup{\cup\cdots\cup}
\newcommand{\vare}{\varepsilon}

\title{Pure pairs. V. Excluding some long subdivision}
\author{Alex Scott\thanks{Research supported by EPSRC grant
EP/V007327/1.}\\
Mathematical Institute, University of Oxford, Oxford OX2 6GG, UK
\\
\\
Paul Seymour\thanks{Supported by AFOSR grant A9550-19-1-0187 and NSF
grant DMS-1800053.}\\
Princeton University, Princeton, NJ 08544
\\
\\
Sophie Spirkl\thanks{This material is based upon work supported by the National Science
Foundation under Award No. DMS-1802201. We acknowledge the support of the Natural Sciences and Engineering Research Council 
of Canada (NSERC) (funding reference
number RGPIN-2020-03912). Cette recherche a \'{e}t\'{e} financée par le Conseil de recherches en sciences naturelles et en 
g\'{e}nie
du Canada (CRSNG) (num\'{e}ro de r\'{e}f\'{e}rence RGPIN-2020-03912).}\\
Princeton University, Princeton, NJ 08544}

\date{December 13, 2019; revised September 27, 2022}

\newtheorem{thm}{}[section]

\newcommand{\Proof}{\noindent{\bf Proof.}\ \ }

\begin{document}
\maketitle
\begin{abstract}
A ``pure pair'' in a graph $G$ is a pair $A,B$ of disjoint subsets of $V(G)$ such that $A$ is complete or anticomplete to $B$.
Jacob Fox showed that for all $\vare>0$, there is a comparability graph $G$ with $n$ vertices,
where $n$ is large,
in which there is no pure pair $A,B$ with $|A|,|B|\ge \vare n$.
He also proved that for all $c>0$ there exists $\vare>0$
such that for every comparability graph $G$ with $n>1$ vertices, there is a pure pair $A,B$
with $|A|,|B|\ge \vare n^{1-c}$; and conjectured that the same holds for every perfect graph $G$.
We prove this conjecture and strengthen it in several ways.

In particular, we show that for all $c>0$, and all $\ell_1, \ell_2\ge 4/c+9$, there exists $\vare>0$ such that, if $G$ is an 
$(n>1)$-vertex graph with no 
hole of length exactly $\ell_1$ and no antihole of length exactly $\ell_2$, then there is a pure pair $A,B$ in $G$
with $|A|\ge \vare n$ and $|B|\ge \vare n^{1-c}$. This is further strengthened, replacing excluding a hole
by excluding some ``long'' subdivision of a general graph.
\end{abstract}

\section{Introduction}
Graphs in this paper are finite, and without loops or parallel edges. Let $A,B\subseteq V(G)$ be disjoint.
We say that $A$ is {\em complete} to $B$, or $A,B$ are {\em complete}, if every vertex in $A$ is adjacent to every vertex in $B$,
and similarly $A,B$ are {\em anticomplete} if no vertex in $A$ has a neighbour in $B$. 
A {\em pure pair} in $G$ is a pair $A,B$ of disjoint subsets of $V(G)$ such that $A,B$
are complete or anticomplete, and $|G|$ denotes the number of vertices of a graph $G$.

Let $\mathcal{H}$ be a set of graphs: we say $G$ is {\em $\mathcal{H}$-free} if no induced subgraph of $G$ is isomorphic to a member 
of $\mathcal{H}$. For some choices of $\mathcal{H}$, every $\mathcal{H}$-free graph admits a pure pair $A,B$ with both $|A|,|B|$
large in terms of $|G|$.
Pure pairs with both $|A|,|B|$ linear in $|G|$ are particularly of interest because of connections with the Erd\H{o}s-Hajnal 
conjecture~\cite{EH0,EH}, and the following was shown in~\cite{pure1}:
\begin{thm}\label{pure1}
Let $\mathcal{H}$ be a finite set of graphs. 
\begin{itemize}
\item If $\mathcal{H}$ contains a forest and the complement of a forest, then there exists $\vare>0$ such that
every $\mathcal{H}$-free graph $G$ with $|G|>1$
admits a pure pair $A,B$ with $|A|,|B|\ge \vare |G|$;
\item If $\mathcal{H}$ does not contain both a forest and the complement of a forest, then there exists $c>0$ and arbitrarily 
large $\mathcal{H}$-free graphs $G$ in which there is no pure pair $A,B$ with $|A|,|B|\ge |G|^{1-c}$.
\end{itemize}
\end{thm}
But if we allow $\mathcal{H}$ to be infinite, the pretty dichotomy of \ref{pure1} disappears: the first bullet remains true, but the second 
may be false. For example, it was shown in~\cite{pure2} that:
\begin{thm}\label{pure2}
Let $H$ be a graph and let $\mathcal{H}$ be the class of all subdivisions of $H$ and their complements; then there exists $\vare>0$ such that
every $\mathcal{H}$-free graph $G$ with $|G|>1$
admits a pure pair $A,B$ with $|A|,|B|\ge \vare |G|$.
\end{thm}
And also, there are classes that do not admit linear pure pairs, but for all $c>0$, do admit pure pairs $A,B$ with $|A|,|B|>|G|^{1-c}$.
For instance, 
Jacob Fox~\cite{fox} proved:
\begin{thm}\label{foxthm}
For every sufficiently large positive integer $n$: 
\begin{itemize}
\item for every $n$-vertex comparability graph $G$,
there is a pure pair $A,B$ in $G$ with
$|A|,|B| > \frac{n}{4\log_2 n}$;
\item there is an $n$-vertex comparability graph $G$
such that
there is no pure pair $A,B$ in $G$ with $|A|,|B|\ge \frac{15n}{\log_2 n}$.
\end{itemize}
\end{thm}
There is also a related asymmetric result, by Fox, Pach and Toth~\cite{toth}:
\begin{thm}\label{toth}
There exists $\vare>0$
such that for every comparability graph $G$ with $|G|>1$, either there is a complete pair $A,B$ with $|A|,|B|\ge \vare |G|$, or
there is an anticomplete pair $A,B$ with $|A|,|B|\ge \vare |G|/\log |G|$.
\end{thm}

Comparability graphs are perfect, and Fox~\cite{fox} (and see also~\cite{foxpach})
conjectured that something like \ref{foxthm} holds for all perfect graphs; more exactly:
\begin{thm}\label{foxconj}
{\bf Conjecture:} For every sufficiently large positive integer $n$ and every $n$-vertex perfect graph $G$,
there is a pure pair $A,B$ in $G$ with
$|A|,|B|\ge n^{1-o(1)}$.
\end{thm}
We will prove this conjecture, and several strengthenings. 
To prove \ref{foxconj} itself, we will show that
\begin{thm}\label{mainthm1}
For all $c>0$, and all sufficiently large $n$, if $G$ is an $n$-vertex perfect graph, then there is
a pure pair $A,B$ in $G$ with $|A|,|B|\ge n^{1-c}$.
\end{thm}
This can be strengthened: we can make one of the two sets of linear size (and replace the 
``sufficiently large'' condition in \ref{mainthm1} with a multiplicative constant). We will show:
\begin{thm}\label{mainthm3}
For all $c>0$ there exists $\vare>0$ such that if $G$ is a perfect graph with $|G|>1$, then there is
a pure pair $A,B$ in $G$ with $|A|\ge \vare|G|$ and $|B|\ge \vare|G|^{1-c}$.
\end{thm}
The complement graph of $G$ is denoted by $\overline{G}$.
A {\em hole} in $G$ is an induced cycle of length at least four, and an {\em antihole} in $G$ is an induced subgraph
whose complement graph is a hole in $\overline{G}$. Perfect graphs are the graphs that have have no holes or antiholes of 
odd length~\cite{SPGT},
but we will show that it is not necessary to exclude all odd holes and odd antiholes to have the result \ref{mainthm3};
it is enough to exclude one of each, of sufficient length. The next result is a strengthening of \ref{mainthm3}:
\begin{thm}\label{mainthm4}
Let $c>0$ with $1/c$ an integer, and let $\ell_1,\ell_2\ge 4/c+5$ be integers. Then there exists $\vare>0$ such 
that 
if $G$ is a graph with $|G|>1$, with no hole of length exactly $\ell_1$ and no antihole of length exactly $\ell_2$, then there is
a pure pair $A,B$ in $G$ with $|A|\ge \vare|G|$ and $|B|\ge \vare|G|^{1-c}$.
\end{thm}

This can be further strengthened, as follows. Let us say $G$ {\em contains} $H$ if some induced subgraph of $G$ is isomorphic to $H$,
and $G$ is  {\em $H$-free} otherwise.
If $X\subseteq V(G)$, $G[X]$ denotes the subgraph induced on $X$.
We say that a graph $H$ has {\em branch-length} at least $\ell$ if every cycle of $H$
has length at least $\ell$, and every two vertices of $H$ with degree at least three have distance at least $\ell$ in $H$.
Since a cycle of length $\ell$ has branch-length $\ell$, the next result strengthens \ref{mainthm4} and is the 
main result of the paper:
\begin{thm}\label{mainthm}
Let $c>0$ with $1/c$ an integer, and let $H_1,H_2$ be graphs with branch-length at least $4/c+5$. Then there exists $\vare>0$ such that 
if $G$ is a graph with $|G|>1$ that is $H_1$-free and $\overline{H_2}$-free, then there is
a pure pair $A,B$ in $G$ with $|A|\ge \vare|G|$ and $|B|\ge \vare|G|^{1-c}$.
\end{thm}

\section{Reduction to the sparse case}

Let us say a graph $G$ is {\em $\vare$-sparse} if every vertex has degree less than $\vare|G|$. We say
$G$ is {\em $(\alpha,\beta)$-coherent} if there do not exist disjoint subsets $A,B$ of $V(G)$, anticomplete to each other,
such that $|A|\ge \alpha$ and $|B|\ge \beta$.

A one-vertex graph does not admit any non-trivial pure pair, but it is $\vare$-sparse for all $\vare>0$,
and $(\alpha,\beta)$-coherent for all $\alpha,\beta>0$; so our standard hypothesis that $G$ is suitably sparse and 
suitably coherent
does not exclude the case $|G|=1$, and we always need to assume separately that $|G|>1$. But we observe:
\begin{thm}\label{big}
Let $0< \vare\le 1/2$; if $G$ is $\vare$-sparse and $(\vare|G|, \vare|G|)$-coherent, with $|G|>1$,
then $|G|> 1/\vare$.
\end{thm}
\Proof
Suppose that $|G|\le 1/\vare$.  If some distinct $u,v\in V(G)$ are non-adjacent, 
$\{u\},\{v\}$ form an anticomplete pair, both of cardinality at least $\vare|G|$, a contradiction.
So $G$ is a complete graph; but its maximum degree is less than $\vare|G|$ and $\vare\le 1/2$, which is impossible since $|G|>1$.
This proves \ref{big}.~\bbox

\bigskip

If $G$ is a graph and $v\in V(G)$, a {\em $G$-neighbour} of $v$ means a vertex of $G$ adjacent to $v$ in $G$.
A theorem of R\"odl~\cite{rodl} implies the following:
\begin{thm}\label{rodlthm}
For every graph $H$ and all $\eta>0$ there exists $\delta>0$ with the following property.
Let $G$ be an $H$-free graph. Then there exists $X\subseteq V(G)$ with $|X|\ge \delta |G|$, such that 
one of $G[X]$, $\overline{G}[X]$ is $\eta$-sparse.
\end{thm}

Consequently, in order to prove \ref{mainthm}, it suffices to prove the following:

\begin{thm}\label{sparsethm}
Let $c>0$ with $1/c$ an integer, and let $H$ be a graph with branch-length at least $4/c+5$. Then there 
exists $\vare>0$ such that 
every $\vare$-sparse $(\vare|G|^{1-c}, \vare|G|)$-coherent graph $G$ with $|G|>1$ contains $H$.
\end{thm}
\noindent{\bf Proof of \ref{mainthm}, assuming \ref{sparsethm}.\ \ }
Let $c>0$ with $1/c$ an integer, and let $H_1,H_2$ have branch-length 
at least $4/c+5$. For $i = 1,2$, choose $\vare_i$ such that \ref{sparsethm} holds with 
$\vare=\vare_i$ and $H=H_i$. Let $\eta=\min(\vare_1,\vare_2,1/2)$. 
Choose $\delta$ such that \ref{rodlthm} holds taking $H=H_1$. Let $\vare=\eta\delta$. We claim that $\vare$ satisfies \ref{mainthm}.

Let $G$ be a graph with $|G|>1$ that is $H_1$-free and $\overline{H_2}$-free.
We must show that there is
a pure pair $A,B$ in $G$ with $|A|\ge \vare|G|$ and $|B|\ge \vare|G|^{1-c}$.
From the choice of $\delta$, there exists $X\subseteq V(G)$ with $|X|\ge \delta |G|$, such that        
one of $G[X]$, $\overline{G}[X]$ is $\eta$-sparse; and by \ref{big} we may assume that $|G|> 1/\vare\ge 1/\delta$, and so $|X|>1$.
In the first case, since $\eta\le \vare_1$, \ref{sparsethm} applied to $G[X]$ 
implies that there is
an anticomplete pair $A,B$ in $G[X]$ with $|A|\ge \eta|X|$ and $|B|\ge \eta|X|^{1-c}$. 
Thus 
$$|A|\ge \eta|X|\ge \eta\delta|G|=\vare|G|$$
and 
$$|B|\ge \eta|X|^{1-c}\ge \eta\delta^{1-c}|G|^{1-c}\ge \eta\delta|G|^{1-c}=\vare|G|^{1-c},$$
as required. In the second case we argue similarly, working in $\overline{G}[X]$. This proves \ref{mainthm}.~\bbox

\bigskip

The remainder of the paper is devoted to proving \ref{sparsethm}.

\section{The pathfinder lemma: finding a path of specified length}

In this section we will prove the main technical tool that we need, which we call the ``pathfinder''.
If $A,B\subseteq V(G)$ are disjoint, we 
say $A$ {\em covers} $B$ if every vertex
in $B$ has a neighbour in $A$. 
A {\em levelling} in $G$ is a sequence $\mathcal{L}=(L_0,L_1\LL L_k)$ of disjoint subsets of $V(G)$ with $k\ge 1$ such that
\begin{itemize}
\item $|L_0|=1$;
\item $L_{i-1}$ covers $L_i$ for $1\le i\le k$; and
\item $L_0\cupcup L_{i-2}$ is anticomplete to $L_i$ for all $i\in \{2\LL k\}$.
\end{itemize}
We denote $L_0\cup L_1\cup\cdots\cup L_k$ by $V(\mathcal{L})$.
We call $L_k$ the {\em base} of the levelling $\mathcal{L}=(L_0,L_1\LL L_k)$, and $V(\mathcal{L})\setminus L_k$ is called the
{\em heart} of $\mathcal{L}$. We call $k$ the {\em height} of the levelling, and the unique vertex in $L_0$ is the {\em apex}.
We call $L_{k-1}$ the {\em penultimate level} of the levelling (for want of a better name).
A path $P$ is {\em $\mathcal{L}$-vertical} if $V(P)\subseteq V(\mathcal{L})$ and $|V(P)\cap L_i|\le 1$ for $0\le i\le k$.

The pathfinder says that if a graph $G$ is suitably sparse and suitably coherent, and we are given
two levellings with disjoint vertex sets and with bases of size linear in $|G|$, and there are suitable constraints on the edges between
the two levellings, then we can find an induced path between the two apexes of any given length greater than the sum of the two heights.
(And there is also a version when the two apexes are equal, and in this case we will find a cycle rather than a path.) 

Let us explain how the pathfinder will be used to prove \ref{sparsethm} and hence \ref{mainthm}. We can assume (by extending $H$ if necessary)
 that the graph $H$ of \ref{sparsethm}
is obtained from some stable set $X$ by adding paths, each of length at least $4/c+5$, where each path has both ends in $X$
and no other vertices in $X$, and all these paths are pairwise vertex-disjoint except for
their vertices in $X$. (For numerical reasons, we will also allow the addition of cycles, but let us skip that for now.)
We are given a graph $G$ which is suitably sparse and suitably coherent, and we need to show that it has an induced subgraph 
isomorphic to $H$.
To obtain a copy of $H$ in $G$, we will choose an appropriate set $X\subseteq V(G)$, and then try to route
paths of $G$ of the right length between the correct pairs of vertices of $X$. We will find each such path by applying the pathfinder to 
some pair of levellings with apexes the corresponding pair of vertices of $X$. But we cannot apply the pathfinder twice to the same levelling,
because the paths we want to produce need to be pairwise disjoint and anticomplete except for their ends. Thus, if some vertex in $X$
is supposed to be an end of several paths of $H$, we will need several levellings all with this apex. So we need a way to
find a good supply of levellings, each with base of linear size, and pairwise disjoint except for their apexes, grouped into
several sets each with a common apex; and we want the edges between them to be under control. And another thing: the pathfinder
can only provide paths between the two apexes of length greater than the sum of the two heights of the corresponding levellings, and we
need paths which might be as short as $4/c+5$, so we need the levellings to have height at most something like $2/c$.

The paper is organized as follows. In this section we prove the pathfinder; and in the next we explain why we can get levellings of 
height about $1/c$ (later, when we try to get several levellings  with a common apex, this height will double). In section 5 we
relax the conditions on levellings, and instead just look for subgraphs of radius about $1/c$ that have a linear set of neighbours (we call this
a ``covering'', and the subgraph of bounded radius is its ``heart''); we find that we can obtain many coverings, with hearts that
are disjoint and pairwise anticomplete. 
What we really want is something slightly more: we want there also to be a vertex
with a neighbour in each of the hearts. To prove this, we prove something stronger, that there is a ``multicovering'', but this is 
just a tool to get one neighbour in common.

So we have many coverings, with hearts pairwise anticomplete and with a common neighbour $a$.
Add $a$ to each of the hearts; then we get many coverings, with hearts pairwise anticomplete except for one common vertex, which we call its ``apex''. We call
this group of coverings a ``spider'', and this is the topic of section 6. 
By making each of the hearts only just big enough that it has linearly many neighbours, we can find a spider
such that most vertices of the graph have no neighbours in any of the hearts of the spider; and so, among them we can do it again,
and find another spider. This way we get a ``troupe'' of spiders, with no edges between their hearts. 
The next step is to convert the hearts of the coverings in each spider to levellings (so the spiders become ``lobsters'');
this is also done in section 6. Then we are ready to apply
the pathfinder, which is done in section 7, and this completes the proof of \ref{sparsethm}.

Let us begin by proving the pathfinder.
First we need the following lemma:

\begin{thm}\label{repeat}
Let $\rho\ge 1$ be some real number, let $K,k>0$ be integers with $K>k$, and let $n_1\LL n_K$ be non-negative integers, 
all less than $\rho^{K/k-2-1/k}$. 
Then there exists $i\in \{1\LL K-k\}$ such that $\rho n_i\ge n_j$ for $j=i+1\LL i+k$.
\end{thm} 
\Proof Suppose not; then for each $i\in \{1\LL K-k\}$ there exists $f(i)$
such that $i<f(i)\le i+k$ and $\rho n_i< n_{f(i)}$. Define $x_1=1$
and $x_{i+1}=f(x_i)$ provided $x_i\le K-k$. Let $x_1\LL x_t$ be defined
by this process; thus $K-k<x_t\le K$. Since $x_{i+1}-x_i\le k$ for each $i$,
it follows that $tk\ge K-1$. Since $n_{x_2}>\rho n_{x_1}$ and $n_{x_2}$
is an integer, it follows
that $n_{x_2}\ge 1$. Thus for $2\le i\le t$, $n_{x_i}\ge \rho^{i-2}$, and so
$n_{x_t}\ge \rho^{t-2}\ge \rho^{K/k-2-1/k}$, contrary to the hypothesis.
This proves \ref{repeat}.~\bbox

Next we need:

\begin{thm}\label{findpath}
Let $c>0$ such that $1/cc$ is an integer, and define $r=2+1/c$.
Let $\ell\ge 1$ be an integer, and define $K=r^{\ell}-1$, and $k=r^{\ell-1}-1$.
Let $\vare>0$, and let $G$ be an $\vare$-sparse $(\vare|G|^{1-c}, \vare|G|)$-coherent graph.
Let $B_0,B_1\LL B_{K}\subseteq V(G)$ be disjoint, where $B_0\ne \emptyset$ and $B_1\LL B_{K}$
each have cardinality
at least $r^{2\ell}\vare|G|$.
Then either:
\begin{itemize}
\item  there is an induced path of length $\ell$,
with vertices $p_0, p_1\LL p_{\ell}$ in order, and 
$$1\le t_1<t_2<\cdots<t_{\ell}\le K,$$
such that $p_0\in B_0$, and $p_i\in B_{t_i}$ for $1\le i\le \ell$; or
\item $|B_0|\le K\vare|G|^{1-c}$, and 
there are sets
$C_1\LL C_{K-k}$ with union $B_0$, such that for 
each $i$ with $1\le i\le K-k$, and each $j$ with $i\le j\le i+k$, at least 
$r^{2\ell-2}\vare|G|$ vertices in $B_j$ have no neighbour in $C_i$.
\end{itemize}
\end{thm}
\Proof
We proceed by induction on $\ell$. Suppose first that $\ell=1$. If there is an edge between $B_0$ and $B_1\cupcup B_K$, then the first bullet holds; 
and if $B_0$ is anticomplete to $B_1\cupcup B_K$,
then since $H$ is $(\vare|G|^{1-c}, \vare|G|)$-coherent and $|B_1|\ge r^{2\ell}\vare|G|\ge \vare|G|$, it follows that $|B_0|<\vare|G|^{1-c}$, and the second bullet
holds, taking $C_1\LL C_{K-k}=B_0$. Thus we may assume that $\ell\ge 2$ and the result holds for $\ell-1$.
Define $\rho=|G|^c$.

Let $B_0=\{v_1\LL v_n\}$.
For all $i\in \{1\LL K-k\}$ and all $j\in \{i\LL i+k\}$, define $A^0_{i,j}=\emptyset$, and $A^0=\emptyset$.
Inductively for $h = 1\LL n$ we will define
\begin{itemize}
\item a set $X^h_i\subseteq B_i$ for each $i\in \{1\LL K\}$
\item the {\em type} of $v_h$ (one of the numbers $1\LL K-k$);
\item a set $A^h_{i,j}\subseteq B_{j}$
for each $i\in \{1\LL K-k\}$ and each $j\in \{i\LL i+k\}$; and
\item a set $A^h$, which is the union of $A^h_{i,j}$ over all $i\in \{1\LL K-k\}$ and all $j\in \{i\LL i+k\}$
\end{itemize}
as follows. Suppose that $1\le h\le n$, and $A^{h-1}$ and $A^{h-1}_{i,j}$ are defined for all $i,j$ with 
$1\le i\le K-k$ and $i\le j\le i+k$. 
For $1\le i\le K$ let $X^{h}_i$ 
be the set
of vertices in $B_i\setminus A^{h-1}$ adjacent to $v_{h}$. 
Since $(K+1)/(k+1)=2+1/c$, it follows that $K>(2+1/c)k+1$, and so $1/c<K/k-2-1/k$.
Hence for $1\le i\le K$, 
$|X^h_i|\le |G|< \rho^{K/k-2-1/k}$. 
By \ref{repeat} applied to the numbers $|X^{h}_1|\LL |X^{h}_K|$,
there exists $t$ with $1\le t\le K-k$ such that $\rho |X^{h}_t|\ge |X^{h}_{j}|$ for $j=t\LL t+k$. Choose some such $t$, which we call
the type of $v_{h}$. 
For each $j\in \{t\LL t+k\}$ define
$A^{h}_{t,j}=A^{h-1}_{i,j}\cup X^{h}_j$; and 
For each $i\in \{1\LL K-k\}\setminus \{t\}$ and each $j\in \{i\LL i+k\}$ define
$A^{h}_{i,j}=A^{h-1}_{i,j}$.
This completes the inductive definition.
\\
\\
(1) {\em  $\rho |A^h_{i,j}|\ge |A^h_{i,i}|$ for all $h \in \{1\LL n\}$ and all $i\in \{1\LL K-k\}$ and all 
$j\in \{i\LL i+k\}$.}
\\
\\
$A^h_{i,j}$ is the disjoint union of the sets $X^h_{j}$ for all $h\in \{1\LL n\}$ such that $v_h$ has type $i$;
and $A^h_{i,i}$ is the disjoint union of $X^h_{i}$ for the same values of $h$. Since $\rho |X^{h}_i|\ge |X^{h}_{j}|$
for each such $h$, this proves (1).
\\
\\
(2) {\em If $v_h$ has type $i$, then every vertex of $B_j$ adjacent to $v_h$
belongs to $A^h$, for all $h\in \{1\LL n\}$, all $1\le i\le K-k$, and all $j\in \{i\LL i+k\}$.}
\\
\\
Let $x\in B_j$ be adjacent to $v_h$.
If $x\in A^{h-1}$, then the claim holds since $A^{h-1}\subseteq A^h$.
If $x\notin A^{h-1}$ then $x\in X^h_j$ from the definition of $X^h_j$, and since $v_h$ has type $i$, it follows that
$$x\in X^h_j\subseteq A^{h}_{i,j}\subseteq A^h.$$
This proves (2).

\bigskip

For $1\le i\le K-k$, let $C_i$ be the set of vertices in $B_0$ that have type $i$. Thus $C_1\LL C_{K-k}$
are pairwise disjoint and have union $B_0$. We note that
$$r^{2\ell}-r^{2\ell-2}=(3+4/c+1/c^2)(k+1)^2\ge 2(k+1)^2.$$
\\
\\
(3) {\em We may assume that $|A^n_{i,j}|> k\vare|G|$ for some $i\in \{1\LL K-k\}$ and some $j\in \{i\LL i+k\}$.}
\\
\\
Suppose not. 
Let $1\le j\le K$. Since $A^{n}\cap B_j$ is the union of the sets $A^n_{i,j}$ for all $i\in \{1\LL K\}$ 
with $j-i\in \{0\LL k\}$, it follows that $|A^n\cap B_j|\le k(k+1)\vare|G|$. 
Now let $1\le i\le K-k$; by (2),
$C_i$ is anticomplete to $B_j\setminus A^n$, for all $j\in \{i\LL i+k\}$. Since
$$|B_j\setminus A^n|=|B_j|-|B_j\cap A^n|\ge r^{2\ell}\vare|G|-k(k+1)\vare|G|\ge r^{2\ell-2}\vare|G|\ge \vare|G|$$
and $G$ is $(\vare|G|^{1-c}, \vare|G|)$-coherent, it follows that $|C_i|<\vare|G|^{1-c}$. 
Hence $|B_0|\le K\vare|G|^{1-c}$; and so the 
second bullet of the theorem holds. This proves (3).

\bigskip

From (3), we may choose $h\in \{1\LL n\}$ minimum such that $|A^h_{i,j}|> k\vare|G|$ for some $i\in \{1\LL K-k\}$ and some $j\in \{i\LL i+k\}$.
Define $D$ to be the set of all $v_{h'}\in C_i$
with $1\le h'\le h$. From the minimality of $h$, and since $G$ is $\vare$-sparse, it follows that 
$|A^h_{i,j}|\le (k+1)\vare|G|$ for all $i\in \{1\LL K-k\}$ and all $j\in \{i\LL i+k\}$.
Consequently $|A^h\cap B_i|\le (k+1)^2\vare|G|$ for $1\le i\le K$.
Now choose $i\in \{1\LL K-k\}$ such that $|A^h_{i,j}|> k\vare|G|$ for some $j\in \{i\LL i+k\}$. By (1), 
$|A^h_{i,i}|> k\vare|G|/\rho= k\vare|G|^{1-c}$.
For $j=i+1\LL i+k$, let $D_j$ be the set of all 
vertices in $B_j$ that have no neighbour in $D$. Thus $B_j\setminus A^h\subseteq D_j$, and so 
$$|D_j|\ge r^{2\ell}\vare|G|-(k+1)^2\vare|G|\ge r^{2\ell-2}\vare|G|.$$
Since $|A^h_{i,i}|> k\vare|G|^{1-c}$, it follows from the inductive hypothesis (with $\ell$ replaced by $\ell-1$, and $B_0$ replaced by $A^h_{i,i}$, and $B_1\LL B_K$
replaced by $D_{i+1}\LL D_{i+k}$) 
that there is an induced path of length $\ell-1$,
with vertices $p_1\LL p_{\ell}$ in order, and 
$$i+1\le t_2<\cdots<t_{\ell}\le i+k,$$
such that $p_1\in A^h_{i,i}$, and $p_i\in B_{t_i}$ for $2\le i\le \ell$. Choose $p_0\in D$ adjacent to $p_1$,
and define $t_1=i$; then the path with vertices $p_0, p_1\LL p_{\ell}$ is induced and satisfies the first bullet of 
the theorem. This proves \ref{findpath}.~\bbox

\bigskip

The ``pathfinder'', the main result of this section, is the following:
\begin{thm}\label{getpath}
Let $c>0$, such that $1/c$ is an integer.
Let $\ell, s,t>0$ be integers, and let $d>0$.
Let $\vare>0$ with $(2+1/c)^{(t+1)(\ell+t)} \vare <d$.
Let $G$ be an $\vare$-sparse $(\vare|G|^{1-c}, \vare|G|)$-coherent graph, and for 
$i = 1,2$, let $\mathcal{L}_i$ be a levelling in $G$, with heart $H_i$,
apex $a_i$, and base $B_i$, satisfying:
\begin{itemize}
\item $V(\mathcal{L}_1)\cap V(\mathcal{L}_2)=\{a_1\}\cap \{a_2\}$; 
\item $V(\mathcal{L}_2)\setminus \{a_2\}$ is anticomplete to $H_1$, 
and if $a_1\ne a_2$ then $a_2$ is anticomplete to $H_1$; 
\item $\mathcal{L}_1, \mathcal{L}_2$ have heights $s,t$ respectively; and
\item  $|B_i|\ge d|G|$ for $i = 1,2$.
\end{itemize}
Then there is an induced path (or cycle, if $a_1=a_2$) of length $\ell+s+t$ between $a_1, a_2$, with vertex set a subset of $H_1\cup B_1\cup H_2\cup B_2$.
\end{thm}
\Proof 
For each integer $i\ge 0$, let $k_{i}=(2+1/c)^{i}-1$.
It follows that 
$k_{\ell}k_{\ell+1}\dots k_{\ell+t}\vare<d$.
Moreover, since $\vare<d$, it follows that $t\ge 2$ (because $|B_2|>\vare|G|$ and $G$ is $\vare$-sparse),
and so $k_{\ell+t} (k_{2\ell+2t}+2)\vare\le d$. 
Define
$d_{i}=(2+1/c)^{2i}\vare$ for each integer $i\ge 0$.

Let $G$, $\mathcal{L}_1$, $\mathcal{L}_2$ and $H_i, a_i, B_i\;(i=1,2)$ satisfy the
hypotheses of the theorem.
Let $\mathcal{L}_1=(L_0\LL L_{s})$ and $\mathcal{L}_2= (M_0\LL M_t)$ say; thus $L_s=B_1$ and $M_t=B_2$.
Let $Z_0=\emptyset$. For $i = 1\LL k_{\ell+t}$, we will inductively define $Z_i\subseteq L_{s-1}$ with $Z_{i-1}\subseteq Z_i$,
and $D_i\subseteq L_s$ with $D_1\LL D_i$ pairwise disjoint, satisfying 
\begin{itemize}
\item $d_{\ell+t}|G|\le |D_i|\le (d_{\ell+t}+\vare)|G|$
\item $D_i$ is the set of all vertices in $L_s$ that have a neighbour in $Z_i$ and have no neighbour in $Z_{i-1}$ 
(and so $D_1\cupcup D_i$ is the set of all vertices in $L_s$ that have a neighbour in $Z_i$).
\end{itemize}
Thus, suppose that $1\le i<k_{\ell+t}$, and 
$Z_0\LL Z_{i-1}$ and $D_1\LL D_{i-1}$ are defined satisfying the conditions above. It follows that
$$|D_1\cupcup D_{i-1}|\le (i-1) (d_{\ell+t}+\vare)|G|\le k_{\ell+t}(d_{\ell+t}+\vare)|G|-d_{\ell+t}|G|.$$
But $d_{\ell+t}+\vare=(k_{2(\ell+t)}+2)\vare$, and $k_{\ell+t} (k_{2\ell+2t}+2)\vare\le d$, so
$$|D_1\cupcup D_{i-1}|\le  k_{\ell+t}(k_{2(\ell+t)}+2)\vare|G|-d_{\ell+t}|G|\le (d-d_{\ell+t})|G|.$$
Hence at least $d_{\ell+t}|G|$ vertices in $L_s$ do not belong to $D_1\cupcup D_{i-1}$. All these vertices have a neighbour
in $L_{s-1}\setminus Z_{i-1}$ and have no neighbour in $Z_{i-1}$; and so there exists $Z_i$ with 
$Z_{i-1}\subseteq Z_i\subseteq L_{s-1}$, minimal such that at least $d_{\ell+t}|G|$ vertices in $L_s$
have a neighbour in $Z_i$ and have none in $Z_{i-1}$. Let this set of vertices be $D_i$. Since $G$ is $\vare$-sparse,
the minimality of $Z_i$ implies that $|D_i|\le (d_{\ell+t}+\vare)|G|$. This completes the inductive definition.

We will try to construct a path (or cycle) satisfying the theorem that starts from $a_2$, runs down through layers of $\mathcal{L}_2$,
jumps to some $D_i$, runs through some of $D_{i+1}, D_{i+2}\LL $ to make it the right length, and then runs up to $a_1$ through
the layers of $\mathcal{L}_1$. The sets $Z_i$ are designed so that when the path has run through enough $D_i$'s
to make its length correct, we can exit into the heart of $\mathcal{L}_1$ without picking up unwanted chords. Note that the only edges 
between $V(\mathcal{L}_2)$ and $V(\mathcal{L}_1)$ have an end in the base of $\mathcal{L}_1$ (or are incident with $a_1$, if $a_1=a_2$).

Let $\mathcal{Q}=(Q_0\LL Q_t)$ be a levelling in $G$. We say it is a {\em sub-levelling} of $\mathcal{L}_2$
if $Q_i\subseteq M_i$ for $0\le i\le t$. For $0\le h\le t$, we say that such a sub-levelling $\mathcal{Q}=(Q_0\LL Q_t)$
is {\em $h$-good} if 
\begin{itemize}
\item there exists $g\in \{1\LL k_{\ell+t}-k_{\ell+t-h}+1\}$, and for each $j\in \{g\LL g+k_{\ell+t-h}-1\}$
there exists $F_j\subseteq D_j$, such that $F_j$ is anticomplete to $Q_0\cup Q_1\cupcup Q_{h-1}$, and $|F_j|\ge d_{\ell+t-h}|G|$; and
\item $|Q_t|>k_{\ell}k_{\ell+1}\dots k_{\ell+t-h}\vare|G|^{1-c}$.
\end{itemize}

Since $d|G|>k_{\ell}k_{\ell+1}\dots k_{\ell+t}\vare|G|^{1-c}$
it follows that $\mathcal{L}_2$ is $0$-good.
Choose $h\le t$ maximum such that some sub-levelling $\mathcal{Q}=(Q_0\LL Q_t\}$ of $\mathcal{L}_2$ is $h$-good,
and let $g$ and the sets $F_j\;(j\in \{g\LL g+k_{\ell+t-h}-1\})$ be as in the definition. 
Let $K=k_{\ell+t-h}$. Since each 
$|F_j|\ge d_{\ell+t-h}|G|$, we may apply \ref{findpath}, replacing $B_0$ by $Q_h$, and replacing $\ell$ by $\ell+t-h$,
and replacing the sequence $B_1\LL B_{k_{\ell}}$ by $F_g\LL F_{g+K-1}$. There are two possible outcomes of \ref{findpath}.

The first outcome is: there is an induced path $P$ of length $\ell+t-h$,
with vertices $p_0, p_1\LL p_{\ell+t-h}$ in order, and
$$g\le t_1<t_2<\cdots<t_{\ell+t-h}\le g+K-1,$$
such that $p_0\in Q_h$, and $p_i\in F_{t_i}$ for $1\le i\le \ell+t-h$. In this case, choose a $\mathcal{Q}$-vertical path $Q$ between 
$a_2$ and $p_0$ (therefore of length $h$); choose a neighbour
$v$ of $p_{\ell+t-h}$ in $Z_{t_{\ell+t-h}}$; and choose an $\mathcal{L}_1$-vertical path $R$ between $a_1, v$ (therefore
of length $s-1$). We claim that 
$$a_2\DD Q\DD p_0\DD P\DD p_{\ell+t-h}\DD v\DD R\DD a_1$$
is an induced path or cycle. To show this, we must check that
\begin{itemize}
\item $V(P)\cap V(Q)=\{p_0\}$, and $V(P)\setminus \{p_0\}$ is anticomplete to $V(Q)\setminus \{p_0\}$; this is true since
$Q_0\LL Q_{h-1}$ are anticomplete to $F_g\LL F_{g+K-1}$ from the definition of $h$-good.
\item $V(P)\cap V(R)=\emptyset$, and the edge with ends $p_{\ell+t-h}$ and $v$ is the only edge between $V(P)$ and $V(R)$;
this is true since $L_0\LL L_{s-2}$ are anticomplete to $L_s$, and $v\in Z_{t_{\ell+t-h}}$ is anticomplete to 
$D_{t_1}\LL D_{t_{\ell+t-h-1}}$.
\item $V(Q)\cap V(R)=\{a_1\}\cap \{a_2\}$, and every edge between $V(Q)$ and $V(R)$ has an end in $\{a_1\}\cap \{a_2\}$; 
this is true from the hypothesis.
\end{itemize}
This proves the path or cycle is indeed induced, and since it has length $\ell$, the theorem holds.

The second outcome of \ref{findpath} is: $\ell+t-h>0$, and $|Q_h|\le K\vare|G|^{1-c}$, and, 
writing $k=k_{\ell+t-h-1}$,
there are sets
$C_g\LL C_{g+K-k-1}$ with union $Q_h$, such that for
each $i$ with $g\le i\le g+K-k-1$, and each $j$ with $i\le j\le i+k$, at least
$d_{\ell+t-h-1}|G|$ vertices in $F_j$ have no neighbour in $C_i$.
Since 
$$|Q_h|\le K\vare|G|^{1-c}< k_{\ell}k_{\ell+1}\dots k_{\ell+t-h}\vare|G|^{1-c}<|Q_t|$$
it follows that $h<t$.
For $g\le i\le g+K-k-1$, let $X_i$ be the set of vertices in $Q_t$ that are joined to a vertex in $C_i$
by a $\mathcal{Q}$-vertical path. Since $\mathcal{Q}$
is a levelling and $C_g\LL C_{g+K-k-1}$ have union $Q_h$, it follows that $X_g\LL X_{g+K-k-1}$ have union $Q_t$;
and since 
$|Q_t|>k_{\ell}k_{\ell+1}\dots k_{\ell+t-h}\vare|G|^{1-c}$,
there exists $i$ with $g\le i\le g+K-k-1$ such that 
$$|X_i|\ge |Q_t|/K>k_{\ell}k_{\ell+1}\dots k_{\ell+t-h-1}\vare|G|^{1-c}.$$
For $h\le h'\le t$ let $Q'_{h'}$ be the set of vertices in $Q_{h'}$ that are joined to a vertex in $C_i$
by a $\mathcal{Q}$-vertical path. Thus $Q'_h=C_i$, and 
$$(Q_0\LL Q_{h-1},Q'_h, Q'_{h+1}\LL Q'_t)$$
is an $(h+1)$-good sub-levelling of $\mathcal{L}_2$, a contradiction. This proves \ref{getpath}.~\bbox

\bigskip

The next result is a form of \ref{getpath} with similar hypotheses except that the bases of the two levellings need not be disjoint, and
we weaken slightly the condition about edges between the heart of $\mathcal{L}_1$ and $V(\mathcal{L}_2)$.
\begin{thm}\label{getpath2}
Let $c>0$, such that $1/c$ is an integer.
Let $\ell, s,t>0$ be integers, and let $d>0$.
Let $\vare>0$ with $(2+1/c)^{(t+1)(\ell+t)} \vare <d/3$.
Let $G$ be an $\vare$-sparse $(\vare|G|^{1-c}, \vare|G|)$-coherent graph.
For $i = 1,2$, let $\mathcal{L}_i$ be a levelling in $G$, with heart $H_i$,
apex $a_i$, and base $B_i$. Suppose that:
\begin{itemize}
\item for $i = 1,2$, $|B_i|\ge d|G|$; 
\item $V(\mathcal{L}_1)\cap V(\mathcal{L}_2)=(\{a_1\}\cap \{a_2\}) \cup (B_1\cap B_2)$; and
\item every edge between $H_1$ and $V(\mathcal{L}_2)$ has one end in the penultimate level of $\mathcal{L}_1$ and the other end in $B_2$.
\end{itemize}
Let $\mathcal{L}_1, \mathcal{L}_2$ have heights $s,t$ respectively.
Then there is an induced path (or cycle, if $a_1=a_2$) of length $\ell+s+t$ between $a_1, a_2$, with vertex set a subset of $H_1\cup B_1\cup H_2\cup B_2$.
\end{thm}
\Proof
Let $d'=d/3$, and 
let $G$, $\mathcal{L}_1$, $\mathcal{L}_2$ and $H_i, a_i, B_i\;(i=1,2)$ satisfy the
hypotheses of the theorem.
Let $\mathcal{L}_1=(L_0\LL L_s)$; thus $L_s=B_1$. 
Choose $L_{s-1}'\subseteq L_{s-1}$ minimal such that at least $d'|G|$ vertices in $B_1\cup B_2$ have a neighbour in $L_{s-1}'$.
Let $L_s'$ be the set of vertices in $B_1\cup B_2$ that have a neighbour in $L_{s-1}'$. Thus 
$$d'|G|\le |L_s'|\le (d'+\vare)|G|\le 2d'|G|.$$

Let $\mathcal{L}_1'$ be the levelling $(L_0\LL L_{s-1}, L_{s-1}', L_s')$. Let $\mathcal{L}_2'$ be the levelling obtained from 
$\mathcal{L}_2$ by replacing its base by $B_2\setminus L_s'$. Then $|L_s'|\ge d'|G|$, and 
$$|B_2\setminus L_s'|\ge d|G|-2d'|G|\ge d'|G|.$$ 
Hence $\mathcal{L}_1'$, $\mathcal{L}_2'$ satisfy the hypotheses of \ref{getpath}, and the result follows. This proves
\ref{getpath2}.~\bbox

\bigskip

When we apply \ref{getpath2}, in the final section, it will be to levellings $\mathcal{L}_1, \mathcal{L}_2$
such that every edge between  $V(\mathcal{L}_1), V(\mathcal{L}_2)$ either is incident with the common apex (when there is one)
or is between the base of one of the levellings and one of the last two terms of the other levelling; so \ref{getpath2} is stronger than we need.

\section{Expansion}

If $X\subseteq V(G)$, $N(X)$ denotes the set of vertices in $V(G)\setminus X$ with a neighbour in $X$, and $N[X]=N(X)\cup X$.
A graph $G$ is {\em $\tau$-expanding} if $|N[X]|\ge \min(\tau|X|,|G|/2)$ for every subset $X\subseteq V(G)$.

\begin{thm}\label{makeexpanding}
Let $c>0$, and let $G$ be a $(|G|^{1-c}/4, |G|/4)$-coherent graph.
Then
there exists $Y\subseteq V(G)$ with $|Y|\le |G|^{1-c}/4$ such that $G\setminus Y$ is $|G|^{c}$-expanding.
\end{thm}
\Proof Let $\alpha=|G|^{1-c}/4$ and $\tau = |G|^{c}$.
Choose $Y\subseteq V(G)$ maximal such that $|Y|\le \alpha$ and $|N[Y]|\le \tau|Y|$ (possibly $Y=\emptyset$).
Let  $W=V(G)\setminus Y$.
If $G[W]$ is $\tau$-expanding then the theorem holds, so we assume not.
Thus
there exists $X\subseteq W$ such that $|N[X]\cap W|< \min(\tau|X|, |W|/2)$. Consequently $X\ne \emptyset$.
But
$$|N[X\cup Y]|\le  |N[Y]|+|N[X]\cap W|\le \tau|Y|+\tau|X|,$$
and so from the maximality of $Y$, it follows that $|X\cup Y|> \alpha$.
Now $|N[Y]|\le \tau|Y|\le \tau\alpha=|G|/4$, and $|N[X]\cap W|\le |W|/2\le |G|/2$; so
$$|N[X\cup Y]|\le |N[Y]| + |N[X]\cap W|\le 3|G|/4.$$
Let $U=V(G)\setminus N[X\cup Y]$; then $|U|\ge |G|/4$.
But $X\cup Y$ is anticomplete to $U$, contradicting that $G$ is $(|G|^{1-c}/4, |G|/4)$-coherent.
This proves \ref{makeexpanding}.~\bbox

\bigskip

If $u,v$ are vertices of a graph $G$, it is sometimes convenient to call the distance between $u,v$ in $G$ the {\em $G$-distance}
between $u,v$.
We deduce:

\begin{thm}\label{smallrad}
Let $c>0$, and let $G$ be a $(|G|^{1-c}/4, |G|/4)$-coherent graph. 
Then there exist $u\in V(G)$ and an integer $k< 1+1/c$, such that:
\begin{itemize}
\item at most $|G|/2$ vertices have $G$-distance less than $k$ from $u$; and
\item at least $|G|/4$ vertices have $G$-distance exactly $k$ from $u$.
\end{itemize}
\end{thm}
\Proof 
By \ref{makeexpanding}, there exists $Y\subseteq V(G)$ with $|Y|\le |G|^{1-c}/4$ 
such that $G\setminus Y$ is $\tau$-expanding, where $\tau = |G|^{c}$. 
Choose $u\in V(G)\setminus Y$, and for each integer $i\ge 0$ let $M_i$ be the set of vertices of $G$
that have $G$-distance at most $i$ from $u$. Since $G\setminus Y$ is $\tau$-expanding,
it follows that for all $i\ge 0$, $|M_{i+1}\setminus Y|\ge \min(\tau|M_i\setminus Y|,|G\setminus Y|/2)$. 
For each $i\ge 1$, let $L_i=M_i\setminus M_{i-1}$.
\\
\\
(1) {\em There exists $k\le 1+1/c$ such that $|L_k|\ge |G|/4$.}
\\
\\
Since $G\setminus Y$ is $\tau$-expanding, it is connected, and so there exists $\ell$ such that $V(G)\setminus Y\subseteq M_\ell$.
Since $|V(G)\setminus Y|\ge 3|G|/4$, we may choose $j\ge 0$ minimum such that $|M_j|\ge |G|^{1-c}/4$. 
Hence for each $i\in \{1\LL j-1\}$, $|M_{i}|< |G|^{1-c}/4<|V(G)\setminus Y|/2$,
and so $|M_{i}\setminus Y| \ge \tau|M_{i-1}\setminus Y|$ since $G\setminus Y$
is $\tau$-expanding. Since $|M_0\setminus Y|=1$,
it follows that $|M_{j-1}\setminus Y|\ge \tau^{j-1}$. Hence
$$|G|^{(j-1)c}=\tau^{j-1}\le |M_{j-1}\setminus Y|\le |M_{j-1}|<|G|^{1-c}/4,$$
and so $(j-1)c < 1-c$, that is, $j< 1/c$.

Since $G$ is $(|G|^{1-c}/4, |G|/4)$-coherent,
and $M_j$ is anticomplete to $V(G)\setminus N[M_j]$, it follows that $|V(G)\setminus N[M_j]|<|G|/4$. But also
$|M_{j-1}|<|G|^{1-c}/4$ (or $j=0$), and so
$|L_j\cup L_{j+1}|\ge |G|-|G|/4 - |G|^{1-c}/4\ge |G|/2$. Thus some $k\in \{j, j+1\}$ satisfies the claim. This proves (1).

\bigskip
Choose $k$ as in (1), minimum. Thus $|L_{k-1}|<|G|/4$, and $|M_{k-2}|< |G|^{1-c}/4$ since 
$G$ is $(|G|^{1-c}/4, |G|/4)$-coherent. Thus $|M_{k-1}|\le |G|/2$.
This proves \ref{smallrad}.~\bbox

\section{Covering sequences}

Let us say a {\em covering} $\mathcal{L}$ in $G$ is a triple $(a,H,B)$ where $H,B$ are disjoint subsets of $V(G)$, $a\in H$, $H$ covers $B$, and 
$G[H]$ is connected. We call $a$ the {\em apex}, $H$ the {\em heart}, and $B$ the {\em base} of the covering, and define 
$V(\mathcal{L})=H\cup B$.
If for every vertex $v\in H$ there is a path of $G[H]$ between $a,v$ of length 
at most $r-1$, we say that $(a,H,B)$ has {\em height} at most $r$, and the least such $r$ is the {\em height} of $(a,H,B)$.
For instance, if $(L_0\LL L_k)$ is a levelling with $k>0$, and $L_0=\{a\}$, then $(a,L_0\cupcup L_{k-1}, L_k)$ is a covering
of height $k$.

A {\em covering sequence} in $G$ is a sequence $(\mathcal{L}_1\LL \mathcal{L}_n)$ of coverings in $G$, with hearts
$H_1\LL H_n$ say, such that $H_1\LL H_n$ are pairwise disjoint and pairwise anticomplete. We call $n$ its {\em length}.
We say such a sequence has {\em height} at most $r$ if each term has height at most $r$. 
If $\mathcal{M}=(\mathcal{L}_1\LL \mathcal{L}_n)$ is a covering sequence, we define $V(\mathcal{M})$ to be the union of the sets
$V(\mathcal{L}_i)$ for $1\le i\le n$.
\begin{figure}[h!]
\centering

\begin{tikzpicture}[scale=0.5,auto=left]

\draw (-5,2) ellipse (2 and 3);
\draw (-0,2) ellipse (2 and 3);
\draw (5,2) ellipse (2 and 3);

\draw (-5,-3) ellipse (3 and 1);
\draw (0,-3) ellipse (3 and 1);
\draw (5,-3) ellipse (3 and 1);

\tikzstyle{every node}=[inner sep=1.5pt, fill=black,circle,draw]
\node (a) at (-5,3) {};
\node (b) at (0,3) {};
\node (c) at (5,3) {};

\draw (-6,1)--(-7,-2.5);
\draw (-5,1)--(-5,-2.5);
\draw (-4,1)--(-3,-2.5);

\draw (-1,1)--(-2,-2.5);
\draw (0,1)--(0,-2.5);
\draw (1,1)--(2,-2.5);

\draw (4,1)--(3,-2.5);
\draw (5,1)--(5,-2.5);
\draw (6,1)--(7,-2.5);

\tikzstyle{every node}=[]
\draw (a) node [above]           {$a_1$};
\draw (b) node [above]           {$a_2$};
\draw (c) node [above]           {$a_3$};
\node at (-5,2) {$H_1$};
\node at (0,2) {$H_2$};
\node at (5,2) {$H_3$};
\node at (-5,-3) {$B_1$};
\node at (0,-3) {$B_2$};
\node at (5,-3) {$B_3$};

\end{tikzpicture}
\caption{A covering sequence of length three. $H_1,H_2,H_3$ are disjoint and anticomplete, but $B_1,B_2,B_3$ need not be;
and there may be edges between $H_i$ and $B_j\setminus B_i$.}
\end{figure}
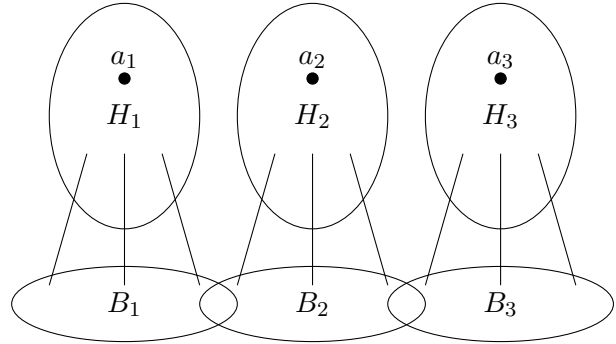

A covering sequence $(\mathcal{L}_1\LL \mathcal{L}_n)$ is a
{\em multicovering} if $\mathcal{L}_1\LL \mathcal{L}_n$ all have the same base, and then this common base is called the {\em base}
of the multicovering.
The main result of this section says that a graph with the usual properties (suitably coherent, suitably sparse) contains a 
multicovering of length any specified constant, with height at most about $1/c$ and with base of linear cardinality. We prove this
in several steps. We begin with:

\begin{thm}\label{getmulticover1}
Let $n\ge 0$ be an integer. Let $c>0$ such that $1/c$ is an integer; let $\vare>0$ with $\vare\le 2^{-n-2}$; and let 
$G$ be an $\vare$-sparse $(\vare|G|^{1-c}, \vare|G|)$-coherent graph.
Then there is a covering sequence $(\mathcal{L}_1\LL \mathcal{L}_n)$ in $G$, where $\mathcal{L}_i=(a_i, H_i, B_i)$ 
for $1\le i\le n$, such that:
\begin{itemize}
\item for $1\le i<j\le n$, $H_i$ is anticomplete to $B_j$;
\item for $1\le i\le n$, $\mathcal{L}_i$ has height at most $1/c$; and
\item for $1\le i\le n$, $|B_i|\ge 2^{-i-1}|G|$.
\end{itemize}
\end{thm}
\Proof We proceed by induction on $n$. If $n=0$ the result is trivial, so we assume that $n\ge 1$ and the result
holds for $n-1$. 
By \ref{smallrad}, there exists $u\in V(G)$ and an integer $k< 1+1/c$ (and hence $k\le 1/c$, 
since $1/c$ is an integer), such that:
\begin{itemize}
\item at most $|G|/2$ vertices of $G$ have distance less than $k$ from $u$; and
\item at least $|G|/4$ vertices of $G$ have distance exactly $k$ from $u$.
\end{itemize}
For $0\le i\le k$ let $L_i$ be the set of all vertices of $G$ with distance exactly $i$ from $u$. Then 
$(L_0\LL L_k)$ is a levelling, with height at most $1/c$; and $|L_k|\ge |G|/4$, so the theorem holds for $n=1$.
Choose $L_{k-1}'\subseteq L_{k-1}$ minimal such that at least $|G|/4$ vertices in $L_k$ have a neighbour in $L_{k-1}'$, and let 
$L_k'$ be the set of vertices in $L_k$ that have a neighbour in $L_{k-1}'$. Thus $|L_k'|\le (1/4+\vare)|G|$ since $G$ is $\vare$-sparse.
Let $\mathcal{L}_1$ be the levelling $(L_0\LL L_{k-2}, L_{k-1}', L_k')$, and let $H_1$ be its heart. Thus 
$|V(\mathcal{L}_1)|\le (3/4+\vare)|G|$.
Let $W$ be the set of vertices of $G$ not in $V(\mathcal{L}_1)$; so $|W|\ge (1/4-\vare)|G|$.
Since $W$ is anticomplete to $H_1$, and $1/4-\vare\ge \vare$
and $G$ is $(\vare|G|^{1-c}, \vare|G|)$-coherent, it follows that $|H_1|\le \vare|G|^{1-c}$, and so 
$|W|\ge (3/4-\vare)|G|-\vare|G|^{1-c}\ge |G|/2$.
Hence $G[W]$ is  $(2\vare)$-sparse and $((2\vare)|W|^{1-c}, (2\vare)|W|)$-coherent.  From the inductive hypothesis applied to $G[W]$,
there is a covering sequence $(\mathcal{L}_2\LL \mathcal{L}_n)$ in $G[W]$, where $\mathcal{L}_i=(a_i, H_i, B_i)$ for $2\le i\le n$,
such that:
\begin{itemize}
\item for $2\le i<j\le n$, $H_i$ is anticomplete to $B_j$;
\item for $2\le i\le n$, $\mathcal{L}_i$ has height at most $1/c$; and
\item for $2\le i\le n$, $|B_i|\ge 2^{-i}|W|\ge 2^{-i-1}|G|$.
\end{itemize}
But then $(\mathcal{L}_1\LL \mathcal{L}_n)$ satisfies the theorem. This proves \ref{getmulticover1}.~\bbox

\bigskip

\begin{thm}\label{getmulticover2}
Let $n\ge 0$ be an integer, let $m=(n-1)^2+1$ and let $\vare= 2^{-m-2}=2^{-n^2+2n-4}$. Let $c>0$ such that $1/c$ is an integer,
and let
$G$ be an $\vare$-sparse $(\vare|G|^{1-c}, \vare|G|)$-coherent graph.
Then there is a covering sequence $(\mathcal{L}_1\LL \mathcal{L}_n)$ in $G$, where $\mathcal{L}_i=(a_i, H_i, B_i)$ for $1\le i\le n$,
such that:
\begin{itemize}
\item for $1\le i\le n$, $\mathcal{L}_i$ has height at most $1/c$;
\item for $1\le i\le n$, $|B_i|\ge 2^{-m-1}|G|$;
\item either $B_1\LL B_n$ are pairwise disjoint and $H_i$ is anticomplete to $B_j$ for all distinct
$i,j\in \{1\LL n\}$, or $B_1=B_2=\cdots=B_n$.
\end{itemize}
\end{thm}
\Proof
Choose $\mathcal{L}_1\LL \mathcal{L}_m$ as in \ref{getmulticover1}, such that each $\mathcal{L}_i$ has base of cardinality
at least $2^{-i-1}|G|$ and height  at most $1/c$. Let $\mathcal{L}_i=(a_i, H_i, B_i)$ for $1\le i\le m$.
For $1\le i\le m$, $H_1\LL H_{i-1}$ are anticomplete to $B_i$, but $H_{i+1}\LL H_m$ might have neighbours in $B_i$.
Choose $B_i'\subseteq B_i$ of cardinality at least $|B_i|/2^{m-i}\ge |G|/2^{m+1}$, such that for each $j\in \{i+1\LL m\}$, 
either every vertex in $B_i'$
has a neighbour in $H_j$, or none do. Let $\mathcal{L}_i'$ be the covering obtained from $\mathcal{L}_i$
by replacing its base by $B_i'$. Then $(\mathcal{L}_1'\LL \mathcal{L}_m')$ is a covering sequence, and for $1\le i<j\le m$,
$H_i$ is anticomplete to $B_j$, and either $H_j$  is anticomplete to $B_i$ or $H_j$ covers $B_i$.
If some $B_i$ is covered by $H_j$ for at least $n$ values of $j$, the theorem holds, so we assume not.
Let $i_1=1$, and inductively for $2\le k\le n$, let $i_k\in \{1\LL m\}$ be minimum such that $H_{i_k}$ is anticomplete
to each of $B_{i_1}'\LL B_{i_{k-1}}'$. This is possible since each of $B_{i_1}'\LL B_{i_{k-1}}'$
is covered by $H_j$ for at most $n-1$ values of $j$, and $m>(n-1)(k-1)$. 
It follows that $B_{i_1}'\LL B_{i_{n}}'$ are pairwise disjoint. Moreover, $i_1<\cdots <i_n$, and so the second outcome holds.
This proves \ref{getmulticover2}.~\bbox

Now we prove the main result of this section. Its proof is closely related to the proof of the main theorem of~\cite{cats}.

\begin{thm}\label{getmulticover3}
Let $c>0$ such that $1/c$ is an integer, and let $n\ge 0$ be an integer. Let $\vare=2^{-2^{2n}}$.
If $G$ is an $\vare$-sparse $(\vare|G|^{1-c}, \vare|G|)$-coherent graph, then 
$G$ contains a multicovering of length $n$ and height at most $1+1/c$, and with base of cardinality at least  
$3\vare |G|$.
\end{thm}

\Proof Define $q=2^{n}$, $p=(q-1)^2+1$, and $x=2^{-p-1}$.
It follows that $\vare\le x3^{-n}$ and $\vare\le 2^{-p-2}$.
From \ref{getmulticover2} (with $m,n$ replaced by $p,q$), we may assume that 
there is a covering sequence $(\mathcal{L}_1\LL \mathcal{L}_q)$ in $G$, 
such that:
\begin{itemize}
\item $V(\mathcal{L}_1)\LL V(\mathcal{L}_q)$ are pairwise disjoint;
\item for $1\le i\le q$, $\mathcal{L}_i$ has height at most $1/c$; 
\item for $1\le i\le q$, the base of $\mathcal{L}_i$ has cardinality at least $x|G|$; and
\item for all distinct $i,j\in \{1\LL q\}$, every edge between $V(\mathcal{L}_i)$ and $V(\mathcal{L}_j)$ is between the base of 
$V(\mathcal{L}_i)$ and the base of $V(\mathcal{L}_j)$.
\end{itemize}
Let $t, d_1\LL d_t> 0$ be integers, where $d_1\LL d_t\le n$. Let us say a {\em battery} with {\em length} $t$ of {\em type} $(d_1\LL d_t)$ is a 
sequence of $t$ multicoverings 
$(\mathcal{M}_1\LL \mathcal{M}_t)$ in $G$, such that:
\begin{itemize}
\item $V(\mathcal{M}_1)\LL V(\mathcal{M}_t)$ are pairwise disjoint;
\item for $1\le i\le t$, $\mathcal{M}_i$ has length $d_i$, and height at most $1+1/c$, and the first term of $\mathcal{M}_i$
has height at most $1/c$;
\item for $1\le i\le t$, the base of $\mathcal{M}_i$ has cardinality at least $x3^{1-d_i}|G|$;
\item for all distinct $i,j\in \{1\LL q\}$, every edge between $V(\mathcal{M}_i)$ and $V(\mathcal{M}_j)$ is between the base of 
$V(\mathcal{M}_i)$ and the base of $V(\mathcal{M}_j)$.
\end{itemize}
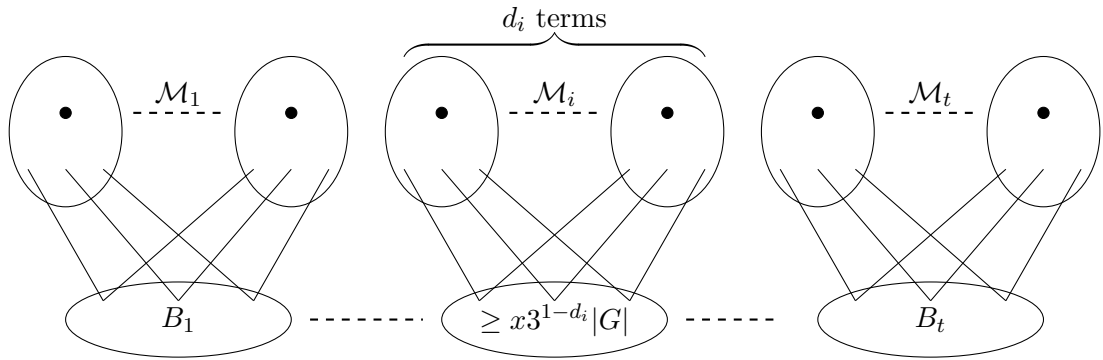
\begin{figure}[h!]
\centering

\begin{tikzpicture}[scale=0.5,auto=left]

\draw (-13,2) ellipse (1.5 and 2);
\draw (13,2) ellipse (1.5 and 2);

\draw (-7,2) ellipse (1.5 and 2);
\draw (7,2) ellipse (1.5 and 2);

\draw (-3,2) ellipse (1.5 and 2);
\draw (3,2) ellipse (1.5 and 2);

\draw (-10,-3) ellipse (3 and 1);
\draw (0,-3) ellipse (3 and 1);
\draw (10,-3) ellipse (3 and 1);

\tikzstyle{every node}=[inner sep=1.5pt, fill=black,circle,draw]
\node (a1) at (-13,2.5) {};
\node (a2) at (-7,2.5) {};
\node (b1) at (-3,2.5) {};
\node (b2) at (3,2.5) {};
\node (c1) at (7,2.5) {};
\node (c1) at (13,2.5) {};

\draw (-14,1)--(-12,-2.5);
\draw (-13,1)--(-10,-2.5);
\draw (-12,1)--(-8,-2.5);

\draw (-8,1)--(-12,-2.5);
\draw (-7,1)--(-10,-2.5);
\draw (-6,1)--(-8,-2.5);

\draw (-4,1)--(-2,-2.5);
\draw (-3,1)--(0,-2.5);
\draw (-2,1)--(2,-2.5);

\draw (2,1)--(-2,-2.5);
\draw (3,1)--(0,-2.5);
\draw (4,1)--(2,-2.5);

\draw (12,1)--(8,-2.5);
\draw (13,1)--(10,-2.5);
\draw (14,1)--(12,-2.5);

\draw (6,1)--(8,-2.5);
\draw (7,1)--(10,-2.5);
\draw (8,1)--(12,-2.5);

\draw[thick, dashed] (-11.2,2.5) -- (-8.8,2.5);
\draw[thick, dashed] (-1.2,2.5) -- (1.2,2.5);
\draw[thick, dashed] (8.8,2.5) -- (11.2,2.5);

\draw[thick, dashed] (-6.5,-3) -- (-3.5,-3);
\draw[thick, dashed] (3.5,-3) -- (6,-3);
\tikzstyle{every node}=[]
\node at (-10,3) {$\mathcal{M}_1$};
\node at (0,3) {$\mathcal{M}_i$};
\node at (10,3) {$\mathcal{M}_t$};
\node at (-10,-3) {$B_1$};
\node at (10,-3) {$B_t$};

\draw [decorate,
    decoration = {calligraphic brace, amplitude = 9pt}, thick] (-4,4) --  (4,4);

\node at (0,5) {$d_i$ terms};
\node at (0,-3) {$\ge x3^{1-d_i}|G|$};
\end{tikzpicture}
\caption{A battery of type $(d_1\LL d_t)$}
\end{figure}

Thus $G$ contains a battery of type $(1\LL 1)$, and of length $q$. Choose a battery $\mathcal{B}$ of type $(d_1\LL d_t)$ with $t$ minimum such that
$2^{d_1}+\cdots +2^{d_t}\ge q$. 
Let $\mathcal{B}=(\mathcal{M}_1\LL \mathcal{M}_t)$. For $1\le i\le t$, let the base of 
$\mathcal{M}_i$ be $B_i$. For each $i$, 
$|B_i|\ge x3^{1-n}|G|\ge 3\vare|G|$.
If some $d_i=n$, then the $i$th term of $\mathcal{B}$ is a multicovering satisfying the theorem; so we assume that $d_1\LL d_t<n$. 
In particular, $2^{d_1}< 2^n=q$, and so $t\ge 2$.
By reordering the 
terms of the battery,
we may assume that $d_t\le d_1\LL d_{t-1}$. 
Since $G$ is $(\vare|G|^{1-c}, \vare|G|)$-coherent, and $|B_t|\ge \vare|G|$, for $1\le i<t$ 
there are fewer than $\vare|G|^{1-c}\le 2|B_i|/3$ vertices in $B_i$ that have no neighbour in $B_t$. Hence
 we may choose $X\subseteq B_t$ minimal such that for some $i\in \{1\LL t-1\}$, at least $|B_i|/3$ vertices in $B_i$
have a neighbour in $X$. For $1\le i<t$, let $Y_i$ be the set of vertices in $B_i$ that have a neighbour in $X$, and $Z_i=B_i\setminus Y_i$.
By reordering, we may assume that $|Y_1|\ge |B_1|/3$. From the minimality of $X$, $|Y_i|\le |B_i|/3+\vare|G|$ for $2\le i\le t-1$, 
and so $|Z_i|\ge 2|B_i|/3-\vare|G|\ge |B_i|/3$.
Let $\mathcal{M}_1=(\mathcal{L}_1\LL \mathcal{L}_{d_1})$, and let the first term of $\mathcal{M}_t$ be $\mathcal{L}=(a,H,B_t)$.
Let $\mathcal{L}_1'$ be the covering $(a, H\cup X, Y_1)$, which therefore has height at most $1+1/c$. Let
$\mathcal{M}_1'$ be obtained from $\mathcal{M}_1$ by replacing its base by $Y_1$ and adding a new final term $\mathcal{L}_1'$; so
$\mathcal{M}_1'$ has length $d_1+1$. For $2\le i\le t-1$, let $\mathcal{M}_i'$ be obtained from $\mathcal{M}_i$ 
by replacing its base by $Z_i$. Then $(\mathcal{M}_1'\LL \mathcal{M}_{t-1}')$ is a battery of type $(d_1+1, d_2\LL d_{t-1})$.
Since $d_1\ge d_t$, it follows that 
$$2^{d_1+1}+\cdots +2^{d_{t-1}}\ge 2^{d_1}+\cdots +2^{d_t}\ge 2^q,$$
a contradiction to the choice of $\mathcal{B}$. This proves \ref{getmulticover3}.~\bbox

\section{Making spiders}

Let $\mathcal{L}_1\LL \mathcal{L}_n$ be coverings in $G$, such that
\begin{itemize}
\item $\mathcal{L}_1\LL \mathcal{L}_n$ all have the same apex $a$;
\item for $1\le i\le n$ let $\mathcal{L}_i$ have heart $H_i$; then for $1\le i<j\le n$, $H_i\setminus \{a\}$ is disjoint from and 
anticomplete to $H_j\setminus \{a\}$.
\end{itemize}
We call $(a,\mathcal{L}_1\LL \mathcal{L}_n)$ a {\em spider} in $G$, and $a$ is its {\em apex}. Its {\em height} is the maximum of the heights of 
$\mathcal{L}_1\LL \mathcal{L}_n$, and its {\em length} is $n$. It has {\em mass} $b$ where $b$ is the minimum cardinality 
of the bases of $\mathcal{L}_1\LL \mathcal{L}_n$. The union of the hearts of $\mathcal{L}_1\LL \mathcal{L}_n$ is called
the {\em heart} of the spider. We call $\mathcal{L}_1\LL \mathcal{L}_n$ the {\em members} of the spider. 
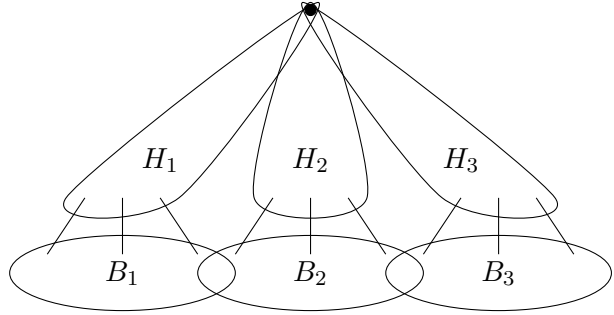
\begin{figure}[h!]
\centering

\begin{tikzpicture}[scale=0.5,auto=left]

\draw plot [smooth cycle] coordinates {(-6.5,-1) (-3.5,-1) (0.2,4.2)};
\draw plot [smooth cycle] coordinates {(-1.5,-1) (1.5,-1) (-0.0,4.2)};
\draw plot [smooth cycle] coordinates {(6.5,-1) (3.5,-1) (-0.2,4.2)};

\draw (-5,-3) ellipse (3 and 1);
\draw (0,-3) ellipse (3 and 1);
\draw (5,-3) ellipse (3 and 1);

\tikzstyle{every node}=[inner sep=1.5pt, fill=black,circle,draw]
\node (a) at (0,4) {};

\draw (-6,-1)--(-7,-2.5);
\draw (-5,-1)--(-5,-2.5);
\draw (-4,-1)--(-3,-2.5);

\draw (-1,-1)--(-2,-2.5);
\draw (0,-1)--(0,-2.5);
\draw (1,-1)--(2,-2.5);

\draw (4,-1)--(3,-2.5);
\draw (5,-1)--(5,-2.5);
\draw (6,-1)--(7,-2.5);

\tikzstyle{every node}=[]
\node at (-4,0) {$H_1$};
\node at (0,0) {$H_2$};
\node at (4,0) {$H_3$};
\node at (-5,-3) {$B_1$};
\node at (0,-3) {$B_2$};
\node at (5,-3) {$B_3$};

\end{tikzpicture}
\caption{A spider of length three.}
\end{figure}

\begin{thm}\label{getspider}
Let $c>0$ such that $1/c$ is an integer, let $n\ge 1$ be an integer, and let $\vare=2^{-2^{2n}}$.
If $G$ is an $\vare$-sparse $(\vare|G|^{1-c}, \vare|G|)$-coherent graph with $|G|\ge 2$, then
$G$ contains a spider of length $n$ and height at most $2+2/c$, and with mass at least
$\vare|G|$.
\end{thm}
\Proof
\ref{getmulticover3} implies that
$G$ contains a multicovering $(\mathcal{L}_1\LL \mathcal{L}_n)$ of length $n$ and height at most $1+1/c$, 
and with base $B$ of cardinality at least $3\vare |G|$. 
Choose $a\in B$. 
Let $1\le i\le n$, and let $H_i$ be the heart of $\mathcal{L}_i$. 
Then every vertex of $H_i\cup \{a\}$ can be joined to $a$
by a path of $G[H_i\cup \{a\}]$ with length at most $1+2/c$. Hence $(a,H_i\cup \{a\}, B\setminus \{a\})$ is a covering of height
at most $2+2/c$, say $\mathcal{L}_i'$.
Consequently $(\{a\}, \mathcal{L}_1'\LL \mathcal{L}_n')$ is a spider of length $n$ and height at most $2+2/c$, and mass
$|B|-1\ge 3\vare |G|-1\ge \vare |G|$, since $|G|\ge 1/\vare$ by \ref{big}. This proves \ref{getspider}.~\bbox

\bigskip

A {\em troupe} of spiders is a set of spiders such that their hearts are pairwise disjoint and anticomplete.

\begin{thm}\label{gettroupe}
Let $c>0$ such that $1/c$ is an integer, and let $m,n \ge 1$ be integers. Let 
let $\vare^{-1} = 2^{2^{2n}}+3(m-1)n$.
If $G$ is an $\vare$-sparse $(\vare|G|^{1-c}, \vare|G|)$-coherent graph, then
$G$ contains a troupe of $m$ spiders, each of length $n$ and height at most $2+2/c$, and with mass at least
$\vare |G|$.
\end{thm}
\Proof
We proceed by induction on $m$.
The result is true if $m=1$, from \ref{getspider}; so we assume that $m\ge 2$, and the result holds for
$m-1$. From \ref{getspider} it follows that 
$G$ contains a spider of length $n$ and height at most $2+2/c$, and with mass at least
$\vare|G|$;  say $\mathcal{S}_1=(a_1,\mathcal{L}_1\LL \mathcal{L}_n)$.
For $1\le i\le n$, let $H_i$ be the heart of $\mathcal{L}_i$. 
Thus every vertex of $H_i$ has $G[H_i]$-distance from $a_1$ at most $1+2/c$, and there are at least $\vare|G|$ vertices in $V(G)\setminus H_i$ with a neighbour in $H_i$.
Let us choose $\mathcal{S}_1$ such that each $H_i$ is
minimal with these two properties (that is, there are at least $\vare|G|$ vertices in $V(G)\setminus H_i$ with a neighbour in $H_i$, 
and every vertex of $H_i$ can be joined to 
$a_1$ by a path of $G[H_i]$ with length at most $1+2/c$.) Let $B_i$ be the set of vertices in $V(G)\setminus H_i$
with a neighbour in $H_i$.
\\
\\
(1) {\em For $1\le i\le n$, $|H_i|< \vare |G|^{1-c}$, and $|B_i|< 2\vare |G|$.}
\\
\\
Let $H_i=\{v_1\LL v_t\}$, ordered with increasing $G[H_i]$-distance from $a_1$ (and hence $v_1=a_1$). Every vertex in $B_i$
either has a neighbour in $H_i\setminus \{v_t\}$ or is adjacent to $v_t$; there are fewer than $\vare|G|$ vertices in $B_i$
with a neighbour in $H_i\setminus \{v_t\}$, from the minimality of $H_i$, and there are fewer than $\vare|G|$ vertices in $B_i$
that are adjacent to $v_t$, since $G$ is $\vare$-sparse. Thus $|B_i|< 2\vare|G|$. Let 
$j=\lceil \vare |G|^{1-c}\rceil$, and suppose that $t\ge j$. Let $J=\{v_1\LL v_j\}$. Thus $|J|\le \vare|G|+1\le 2\vare|G|$ 
by \ref{big}, and so 
$|V(G)\setminus J|\ge (1-2\vare)|G|$. Since $G$ is 
$(\vare|G|^{1-c}, \vare|G|)$-coherent, there are fewer than $\vare |G|$ vertices in $V(G)\setminus J$ that have no neighbour in $J$,
and so there are at least $(1-3\vare)|G|\ge 2\vare|G|$ vertices in $V(G)\setminus J$ that have a neighbour in $J$. 
Thus $t=j$ and $H_i=J$, from the minimality of $H_i$; but this is impossible since $|B_i|<2\vare|G|$. This proves that $t<j$, and so proves (1).

\bigskip

From (1), 
$$|H_1\cupcup H_n\cup B_1\cupcup B_n|\le 3\vare n|G|.$$
Let $X$ be the complement in $V(G)$ of this set; thus $|X|\ge (1-3\vare n)|G|$.
Let 
$$(\vare')^{-1} = 2^{2^{2n}}+3(m-2)n=\vare^{-1} - 3n.$$ 
Thus 
$\vare'(1-3n\vare)=\vare$, and so $\vare'|X|\ge \vare|G|$.
It follows that $G[X]$ is $\vare'$-sparse and $(\vare'|X|^{1-c}, \vare'|X|)$-coherent.
From the inductive hypothesis applied to $G[X]$, we deduce
that there is a troupe of $m-1$ spiders in $G[X]$, each of length $n$ and height at most $2+2/c$, and with mass at least
$\vare'|X|\ge \vare|G|$. But then adding $\mathcal{S}_1$ to this troupe gives a troupe of $m$ spiders satisfying the theorem.
This proves \ref{gettroupe}.~\bbox

\bigskip

So, our graph contains a troupe of spiders, of arbitrarily large cardinality, and each with arbitrarily large length, all of height
at most $2+2/c$, and with bases of linear cardinality.
The next result converts the members of these spiders to levellings, but we need to be careful exactly what we mean. In a levelling,
all edges from heart to base start from the penultimate level of the levelling. We need more than this: we need that for every two 
levellings that are members of spiders in the troupe, every edge from the heart of one to the base of 
the other starts from the penultimate level of the first, and this is more tricky to arrange. 

Let us first state the definition formally. Let $n\ge 1$ and let $\mathcal{L}_1\LL \mathcal{L}_n$ be levellings in a graph 
$G$, all with the same apex $a$,
such that 
\begin{itemize}
\item for $1\le i\le n$, let $H_i$ be the heart of $\mathcal{L}_i$; then $H_1\setminus \{a\}\LL H_n\setminus \{a\}$
are pairwise disjoint (the bases of $\mathcal{L}_1\LL \mathcal{L}_n$ may intersect);
\item for all distinct $i,j\in \{1\LL n\}$, every edge of $G$ between $H_i\setminus \{a\}$ and $V(\mathcal{L}_j)\setminus \{a\}$ 
has one end in the penultimate level of $\mathcal{L}_i$ and the other in the base of $\mathcal{L}_j$.
\end{itemize}
We call $(a,\mathcal{L}_1\LL \mathcal{L}_n)$ a
{\em lobster} in $G$, and $a$ is its {\em apex}. Its {\em height} is the maximum height of $\mathcal{L}_1\LL \mathcal{L}_n$,
and its {\em length} is $n$. It has {\em mass} $b$ where $b$ is the minimum cardinality
of the bases of $\mathcal{L}_1\LL \mathcal{L}_n$. Its {\em heart} is the union of the 
hearts of $\mathcal{L}_1\LL \mathcal{L}_n$. We call $\mathcal{L}_1\LL \mathcal{L}_n$ the {\em members} of the lobster.

\begin{figure}[h!]
\centering

\begin{tikzpicture}[scale=0.5,auto=left]

\draw (-5,-3) ellipse (3 and 1);
\draw (0,-3) ellipse (3 and 1);
\draw (5,-3) ellipse (3 and 1);

\draw (-5,-.5) ellipse (2 and .7);
\draw (-0,0) ellipse (2 and .7);
\draw (5,-1) ellipse (2 and .7);

\draw (-3.5,1.5) ellipse (1.5 and .5);

\draw (4,.5) ellipse (1.5 and .5);
\draw (3,1.7) ellipse (1 and .3);

\tikzstyle{every node}=[inner sep=1.5pt, fill=black,circle,draw]
\node (a) at (0,3) {};

\draw (-6,-.8)--(-7,-2.5);
\draw (-5,-.8)--(-5,-2.5);
\draw (-4,.-.8)--(-3,-2.5);

\draw (-1,-.3)--(-2,-2.5);
\draw (0,-.3)--(0,-2.5);
\draw (1,-.3)--(2,-2.5);

\draw (4,-1.3)--(3,-2.5);
\draw (5,-1.3)--(5,-2.5);
\draw (6,-1.3)--(7,-2.5);

\draw (a) -- (-1,.3);
\draw (a) -- (0,.3);
\draw (a) -- (1,.3);

\draw (a) -- (-2.7,1.8);
\draw (a) -- (-3.5,1.8);
\draw (a) -- (-4.3,1.8);

\draw (-2.7,1.3) -- (-3.5,-.2);
\draw (-3.5,1.3) -- (-4.8,-.2);
\draw (-4.3,1.3) -- (-6.1,-.2);

\draw (a) -- (2.5,1.9);
\draw (a) -- (3,1.9);
\draw (a) -- (3.5,1.9);

\draw (2.5,1.6) -- (3,.7);
\draw (3,1.6) -- (4,.7);
\draw (3.5,1.6) -- (5,.7);

\draw (3,.3)-- (4,-.7);
\draw (4,.3) -- (5,-.7);
\draw (5,.3) -- (6,-.7);

\tikzstyle{every node}=[]
\draw (a) node [above]           {$a$};

\end{tikzpicture}
\caption{A lobster of length three.}
\end{figure}
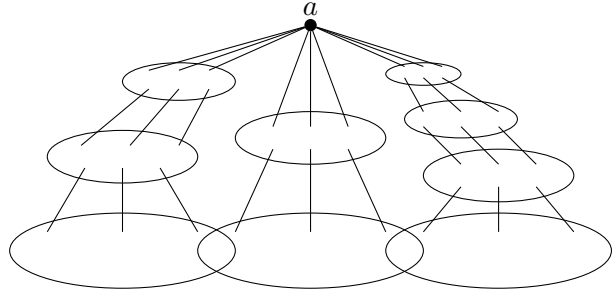

A {\em troupe} of lobsters is a set $\{\mathcal{T}_1\LL \mathcal{T}_m\}$ of lobsters, such that for all $i,j\in \{1\LL m\}$:
\begin{itemize}
\item for  $1\le i<j\le m$, the heart of $\mathcal{T}_i$ is disjoint from and anticomplete to the heart of $\mathcal{T}_j$;
\item let $\mathcal{L}, \mathcal{M}$ each be a member of one of $\mathcal{T}_1\LL \mathcal{T}_m$, with $\mathcal{L}\ne \mathcal{M}$,
and let $\mathcal{L}=(L_0\LL L_k)$; then there is no edge between $L_0\cupcup L_{k-2}$ and the base of $\mathcal{M}$.
\end{itemize}

\begin{thm}\label{spiderlobster}
Let $c>0$ such that $1/c$ is an integer, and let $m,n\ge 1$ be integers. Let 
$$\vare^{-1}= (2^{2^{2n}}+3(m-1)n)(2+2/c)^{mn}.$$
If $G$ is an $\vare$-sparse $(\vare|G|^{1-c}, \vare|G|)$-coherent graph, then
$G$ contains a troupe of $m$ lobsters, each of length $n$ and height at most $2+2/c$, and with mass at least
$\vare|G|$.
\end{thm}
\Proof 
Let $G$ be an $\vare$-sparse $(\vare|G|^{1-c}, \vare|G|)$-coherent graph. By \ref{gettroupe} with $\vare$ replaced by 
$(2+2/c)^{mn}\vare$,
there is a troupe of spiders $\{\mathcal{S}_1\LL \mathcal{S}_m\}$ in $G$, each of length $n$ and height at most $2+2/c$, 
and with mass at least $(2+2/c)^{mn}\vare|G|$. Let the members of these spiders (in some order) be $\mathcal{L}_1\LL \mathcal{L}_{mn}$,
and for $1\le i\le mn$ let $\mathcal{L}_i=(a_i,H_i, B_i)$. (Thus, some of $a_1\LL a_{mn}$ may be equal.)
We shall convert these members one by one to levellings, at each step 
shrinking all the bases.

Let $X^0=B_1\cupcup B_{mn}$, and for $1\le i\le mn$ let $X^0_i$ be the set of all vertices in $X^0$ with a neighbour in $H_i$
(thus $B_i\subseteq X^0_i$). Inductively, let $1\le h\le mn$, and suppose that we have defined $X^{h-1}$ and 
$\mathcal{L}'_1\LL \mathcal{L}_{h-1}'$, and $X^{h-1}_i$ for $1\le i\le mn$, satisfying:
\begin{itemize}
\item for $1\le i\le h-1$, $\mathcal{L}'_i$ is a levelling; its heart is a subset of $H_i$, and $a_i$ is its apex; its height
is at most $2+2/c$;
\item for $1\le i\le h-1$, $X^{h-1}_i$ is the set of all vertices in $X^{h-1}$ with a neighbour in the heart of $\mathcal{L}'_i$,
and for $h\le i\le mn$, $X^{h-1}_i$ is the set of all vertices in $X^{h-1}$ with a neighbour in the heart of $\mathcal{L}_i$;
\item for $1\le i\le h-1$, every edge between the heart of $\mathcal{L}'_i$ and $X^{h-1}$ has an end in the penultimate level of 
$\mathcal{L}'_i$; and 
\item for $1\le i\le mn$, $|X^{h-1}_i|\ge (2+2/c)^{mn+1-h}\vare|G|$.
\end{itemize}

For $0\le j\le 1+2/c$, let $L_j$ be the set of vertices
in $H_h$ with $G[H_h]$-distance to $a_h$ exactly $j$. Thus every vertex $v\in X^{h-1}_h$ has a neighbour in some
$L_j$ where $j\in \{0\LL 1+2/c\}$, and the smallest such $j$ is called the {\em type} of $v$. There are only 
$2+2/c$ possible types, and so there exists $k\in \{0\LL 1+2/c\}$ such that at least $|X^{h-1}_h|/(2+2/c)$
vertices in $X^{h-1}_h$ have type $k$. Consequently, since 
$$|X^{h-1}_h|/(2+2/c)\ge (2+2/c)^{mn+1-h}\vare|G|/(2+2/c)= (2+2/c)^{mn-h}\vare|G|,$$
there exists $k\in  \{0\LL 1+2/c\}$ minimum such that at least $(2+2/c)^{mn-h}\vare|G|$
vertices in $X^{h-1}_h$ have type $k$. Let $X^h_h$ be the set of all vertices in $X^{h-1}_h$ that have type $k$, and 
let $\mathcal{L}'_h=(L_0\LL L_k,X^h_h)$. Thus $\mathcal{L}'_h$ is a levelling with height $k+1\le 2+2/c$.

Let $Z^h$ be the set of vertices in $X^{h-1}_h$ with type less than $k$.
Thus 
$$|Z^h|\le (1+2/c)(2+2/c)^{mn-h}\vare|G|.$$
For $1\le i\le mn$ with $i\ne h$, define $X^h_i=X^{h-1}_i\setminus Z^h$. Thus $|X^h_i|\ge |X^{h-1}_i|-|Z^h|$,
and so 
$$|X^h_i|\ge (2+2/c)^{mn+1-h}\vare|G|- (1+2/c)(2+2/c)^{mn-h}\vare|G|= (2+2/c)^{mn-h}\vare|G|.$$
Let $X^h$ be the union of the sets $X^h_i\;(1\le i\le n)$.
This completes the inductive definition.

For $1\le i\le m$, let $\mathcal{T}_i$ be the lobster obtained from $\mathcal{S}_i$ by replacing each of its members $\mathcal{L}_j$
by $\mathcal{L}_j'$. This makes a troupe of lobsters satisfying the theorem, and so proves \ref{spiderlobster}.~\bbox

\section{Part assembly}
Now we put these several pieces together to prove \ref{mainthm}, which we restate:
\begin{thm}\label{mainthmagain}
Let $c>0$ with $1/c$ an integer, and let $H_1,H_2$ be graphs with branch-length at least $4/c+5$. Then there exists $\vare>0$ such that
if $G$ is a graph with $|G|>1$ that is $H_1$-free and $\overline{H_2}$-free, then there is
a pure pair $A,B$ in $G$ with $|A|\ge \vare|G|$ and $|B|\ge \vare|G|^{1-c}$.
\end{thm}

As we saw in section 2, to prove \ref{mainthmagain}, it suffices to prove the following:

\begin{thm}\label{sparsethmagain}
Let $c>0$ with $1/c$ an integer, and let $H$ be a graph with branch-length at least $4/c+5$. Then there exists $\vare>0$ such that
if $G$ is an $\vare$-sparse $(\vare|G|^{1-c},\vare|G|)$-coherent graph, then $G$ contains $H$.
\end{thm}
\Proof
By adding more vertices to $H$, we may assume that if $X$ denotes the set of vertices of $H$ that have degree different from two,
then every cycle of $H$ contains at least one vertex in $X$, and every path in $H$ with both ends in $X$ has length at least 
$4/c+5$, and every cycle of $H$ has length at least $4/c+5$. Let $X=\{x_1\LL x_m\}$. 
Consequently $H$ can be obtained from the set $X$
of $m$ isolated vertices by adding 
\begin{itemize}
\item paths with ends in $X$ and each of length at least $4/c+5$, and
\item cycles with exactly one vertex in $X$, of length at least $4/c+5$
\end{itemize}
where every vertex of $V(H)\setminus X$ belongs to exactly one of these paths and cycles, and has degree exactly two in $H$.
Let the paths be $R_i\; (i\in I_1)$, and let the cycles be $R_i\;(i\in I_2)$, where $I_1\cap I_2=\emptyset$. For $i\in I_1$,
let $R_i$ have ends $(u_i, v_i)$ (ordered arbitrarily) and have length $\ell_i$, and for $i\in I_2$, 
let $u_i=v_i$ be the unique vertex of $R_i$ in $X$, and let $R_i$ have length $\ell_i$. 
Thus $H$ is determined up to isomorphism
by a knowledge of $X$, the pairs $(u_i,v_i)\;(i\in I_1\cup I_2)$, and the 
numbers $\ell_i\;(i\in I_1\cup I_2)$. Let $I=I_1\cup I_2$, and for each $i\in I$ let $\alpha_i, \beta_i\in \{1\LL m\}$ such that $x_{\alpha_i}=u_i$
and $x_{\beta_i}=v_i$. Let $I=\{1\LL p\}$.

Let $n$ be the maximum degree of $H$, and let  
$d^{-1}= (2^{2^{2n}}+3(m-1)n)(2+2/c)^{mn}$.
Choose $\vare>0$ with $3(2+1/c)^{|H|^2}(4p)^p \vare <d$.
We claim that
$\vare$ satisfies the theorem. Let $G$ be an $\vare$-sparse $(\vare|G|^{1-c}, \vare|G|)$-coherent graph. We must
show that $G$ contains $H$. By \ref{spiderlobster},
$G$ contains a troupe $\{\mathcal{S}_1\LL \mathcal{S}_m\}$ of $m$ lobsters, each of length $n$ and height at 
most $2+2/c$, and with mass at least $d|G|$. 
For $1\le i\le p$, choose a member $\mathcal{L}_{2i-1}$ of $\mathcal{S}_{\alpha_i}$
and a member $\mathcal{L}_{2i}$ of $\mathcal{S}_{\beta_i}$, in such a way that the levellings $\mathcal{L}_i, \mathcal{M}_i\;(i\in I)$ are all different.
(This is possible from the definition of $n$).

We will prove that for all $h\in I$, there is a path $P_h$ (or cycle, if the two apexes are equal) between the
apex of $\mathcal{L}_{2h-1}$ and the apex of $\mathcal{L}_{2h}$ of length $\ell_h$, such that the union of  $P_1\LL P_p$ 
makes an induced subgraph of
$G$ isomorphic to $H$. We will choose these paths and cycles in order. Also for $1\le h\le p$ we need to choose 
a subset $X^h_{k}$ of the base of each $\mathcal{L}_k$ for $2h<k\le 2p$, and a subset $Y^h_k$ of the penultimate level of 
$\mathcal{L}_k$, with properties that we will describe below.
We denote by $P_h^*$ the set of vertices of $P_h$ different from its ends (if it is a path)
or different from the apex of $A_h$ (if it is a cycle). In either case $|P_h^*|=\ell_h-1$.

For $0\le h\le p$ let $w_h=(4p)^{-h}d$.
Let $B$ be the union of the bases of $\mathcal{L}_1\LL\mathcal{L}_{2\rho}$. 
For $1\le i\le 2p$, let 
$Y^0_i$ be the penultimate level of $\mathcal{L}_i$, let $X^0_i$ be the set of vertices in $B$ with a neighbour in $Y^0_i$, and let
$a_i$ be the apex of $\mathcal{L}_i$. Thus, $|X^0_i|\ge w_0|G|$.
Now inductively, suppose we have chosen
the first $h-1$ paths or cycles, say $P_1\LL P_{h-1}$, where $1\le h\le p$, satisfying:
\begin{itemize}
\item for $1\le g\le h-1$, if $a_{2g-1}\ne a_{2g}$, then $P_g$ is an induced path joining these apexes, 
of length $\ell_g$; and if the apexes are equal then $P_g$ is a cycle of length $\ell_g$ containing this apex;
\item for $1\le g\le h-1$, and $2h+1\le i\le 2p$, every vertex of the heart of $\mathcal{L}_{i}$ with a neighbour in $P_g^*$ belongs to
the penultimate level of $\mathcal{L}_{i}$.
\end{itemize}
Suppose moreover that for $2h-1\le i\le 2p$ we have chosen $X^{h-1}_{i}\subseteq X^0_{i}$ and $Y^{h-1}_i\subseteq Y^0_{i}$, such that
for all $i\in \{2h-1\LL 2p\}$:
\begin{itemize}
\item $X^{h-1}_i$ is the set of all vertices in $B$ with a neighbour in $Y^{h-1}_i$;
\item $X^{h-1}_i\cup Y^{h-1}_i$ is anticomplete to $P_1^*\LL P^*_{h-1}$; and
\item $|X^{h-1}_{i}|\ge w_{h-1}|G|$.
\end{itemize}
We choose $P_h$ as follows. For $2h+1\le i\le 2p$, choose $Y_i^h\subseteq Y^{h-1}_i$ minimal such that at least $(w_h+\vare(|H|-1))|G|$ vertices 
in $B$ (necessarily all in $X^{h-1}_i$) have a neighbour in $Y_i$, and let $X_i$ be the set of vertices in $B$ with a neighbour
in $Y_i$. From the minimality of $Y_i$, 
$$(w_h+\vare(|H|-1))|G|\le |X_i|\le (w_h+\vare|H|)|G|.$$ 
Now $\vare|H| \le w_h$, so $|X_i|\le 2w_h$.
Let $Z=X_{2h+1}\cupcup X_{2p}$; 
thus $|Z|\le 2(2p-2)w_h|G|$. For $i = 2h-1, 2h$ let $X_i=X^{h-1}_i\setminus Z$. Thus 
$$|X_i|\ge |X^{h-1}_i|-|Z|\ge (w_{h-1} -4(p-1)w_h)|G|\ge w_h|G|\ge (4p)^{-p}d|G|= 3(2+1/c)^{|H|^2}\vare|G|$$
for $i = 2h-1, 2h$. 

For $i = 2h-1,2h$ let $\mathcal{L}_i'$ be the levelling obtained from $\mathcal{L}_i$
by replacing its base by $X_i$. 
Now $\mathcal{L}_{2h-1}'$, $\mathcal{L}_{2h}'$ both have height at most $2+2/c$, and $\ell_h\ge 5+4/c$.
By \ref{getpath2} applied to the levellings $\mathcal{L}_{2h-1}'$, $\mathcal{L}_{2h}'$, there is an induced path $P_h$
of length $\ell_h$ between $a_{2h-1}, a_{2h}$ (or a cycle, if $a_{2h-1}=a_{2h}$), with vertex set included in 
$V(\mathcal{L}_{2h-1}')\cup V(\mathcal{L}_{2h}')$. Consequently $P_h^*$ is anticomplete to $Y_i$ for $2h+1\le i\le 2p$,
and to $P_1^*\LL P_{h-1}^*$. It might have neighbours in $X_i$ for $2h+1\le i\le 2p$, but since $|P_h^*|\le |H|-1$,
there are at most $\vare(|H|-1)|G|$ such vertices. For $2h+1\le i\le 2p$, let $X^h_i$ be the set of vertices in $X_i$
with no neighbour in $P_h^*$. Thus $|X^h_i|\ge |X_i|-\vare(|H|-1)|G|\ge w_h|G|$. This completes the inductive definition.

But then the union of $P_1\LL P_p$ forms an induced subgraph isomorphic to $H$. This proves \ref{sparsethmagain},
and hence completes the proof of \ref{mainthm}.~\bbox

One might wonder how $\vare$ in \ref{mainthm} depends on $c, H_1,H_2$.
For simplicity let us assume that 
if $X$ denotes the set of vertices of $H_1$ that have degree different from two,
then every cycle of $H_1$ contains at least one vertex in $X$, and every path in $H_1$ with both ends in $X$ has length at least
$4/c+5$, and every cycle of $H_1$ has length at least $4/c+5$; and the same for $H_2$. Let $r=\max(|H_1|,|H_2|)$.
Then one can check that (with $H=H_1,H_2$) defining $d,\vare$
as in the proof of \ref{sparsethmagain} yields a value of $\vare$ that satisfies \ref{sparsethmagain}, with $\log\log(1/\vare)=O(r)$. Next we need a version of
\ref{rodlthm} with an explicit dependence of $\delta$ on $\eta$, and for this we can use the proof of \ref{rodlthm} due to 
Fox and Sudakov~\cite{foxsudakov}; this can be used to show that \ref{rodlthm} holds where $\log 1/\delta=O(|H|(\log 1/\eta)^2)$.
The argument given in section 2 that \ref{sparsethm} implies \ref{mainthm} shows that if $\vare$ satisfies \ref{sparsethmagain}, 
and $\delta$ satisfies \ref{rodlthm}, then 
$\vare'=\vare\delta$ satisfies \ref{mainthm}. Putting these pieces together, we deduce that there exists $\vare$ with 
$\log\log(1/\vare)=O(r)$ that satisfies \ref{mainthm}.
\section{Further extension}

We have found a kind of strengthening of \ref{mainthm}, that we state without proof.
For $\ell\ge 2$, let us say a graph $H$ is {\em $\ell$-handled}
if there are induced subgraphs $P_0\LL P_k$ of $H$, for some $k\ge 1$, such that:
\begin{itemize}
\item $P_0$ is a forest;
\item every path of $P_0$ has length at most $\ell$;
\item $P_1\LL P_k$ are pairwise vertex-disjoint paths, each of length at least $\ell$;
\item for $1\le i\le k$, $V(P_i\cap P_0)$ consists exactly of the two ends of $P_i$; and
\item $H=P_0\cup P_1\cupcup P_k$.
\end{itemize}
Then:
\begin{thm}\label{treethm}
There exists $\gamma>0$ with the following property.
Let $c>0$ with $1/c$ an integer, and let $H_1,H_2$ be $\gamma/c $-handled graphs. Then there exists $\vare>0$ such that
if $G$ is a graph with $|G|>1$ that is $H_1$-free and $\overline{H_2}$-free, then there is
a pure pair $A,B$ in $G$ with $|A|\ge \vare|G|$ and $|B|\ge \vare|G|^{1-c}$.
\end{thm}
Then the essentials of \ref{mainthm} follow from \ref{treethm} by taking $P_0$ to be the subgraph of $H$ induced on the set of
all vertices of degree
at least three and their neighbours. But we feel that \ref{treethm} is not very satisfactory, because if the forest $P_0$
has long paths, the hypothesis requires the paths $P_1\LL P_k$ to be long too. We would prefer a version of \ref{treethm} where we omit the
second bullet from the definition of $\ell$-handled, but so far we cannot prove it.

A weaker form of \ref{mainthm} will be proved for a wider class of graphs in~\cite{pure8}. Let $H$ be a graph.
If $E(H)\ne \emptyset$, we define the {\em congestion}
of $H$ to be the maximum of $1-(|J|-1)/|E(J)|$, taken over all subgraphs
$J$ of $H$ with at least one edge; and if $E(H)=\emptyset$, we define
the congestion of $H$ to be zero. Thus the congestion of $H$ is always non-negative, and equals zero
if and only if $H$ is a forest; and, for instance, long cycles have smaller congestion than short cycles.

In~\cite{pure8} we will prove:
\begin{thm}\label{congthm}
Let $c>0$, and let $H_1,H_2$ be graphs with congestion at most $c/(9+15c)$.
Then there exists $\vare>0$ such that
if $G$ is a graph with $|G|>1$ that is $H_1$-free and $\overline{H_2}$-free, then there is
a pure pair $A,B$ in $G$ with $|A|,|B|\ge \vare|G|^{1-c}$.
\end{thm}

This is pleasing because of the following weak converse (easily proved with a random graph argument that we omit):
\begin{thm}\label{conjconverse}
Let $c>0$, and let $H_1,H_2$ be graphs both with congestion more than $c$.
There is no $\vare>0$ such that
for every graph $G$ with $|G|>1$ that is $H_1$-free and $\overline{H_2}$-free, there is
a pure pair $A,B$ in $G$ with $|A|, |B|\ge \vare|G|^{1-c}$.
\end{thm}

The result \ref{congthm} does not contain \ref{mainthm}, because in \ref{congthm} neither of $A,B$ have to have linear
cardinality. What if we ask for a strengthened version of \ref{congthm} that would contain \ref{mainthm} (by requiring one of $|A|,|B|$ to be linear)?
We pose that as a conjecture:
\begin{thm}\label{conj}
{\bf Conjecture:} For all $c>0$, there exists $\sigma>0$ with the following property. Let $H_1,H_2$ be graphs with congestion at most $\sigma$.
There exists $\vare>0$ such that
if $G$ is a graph with $|G|>1$ that is $H_1$-free and $\overline{H_2}$-free, then there is
a pure pair $A,B$ in $G$ with $|A|\ge \vare|G|$ and $|B|\ge \vare|G|^{1-c}$.
\end{thm}

\section*{Acknowledgement}
The referee reports were very thorough and very hepful, and we would like to express our thanks.

\end{document}